%% file: main.tex
\documentclass[a4paper,12pt,parskip=half]{amsart}

\usepackage[foot]{amsaddr}
\usepackage{psfrag}
\usepackage{subfigure}
\usepackage{minitoc}
\usepackage{tikz}
\usepackage{enumerate}
\usepackage{filecontents}
\usepackage[margin=1 in]{geometry}
\usepackage{amsmath}

\newcommand{\n}{{\bf n}}
\newcommand{\x}{{\bf x}}
\newcommand{\X}{{\bf X}}
\newcommand{\f}{{\bf f}}
\newcommand{\xiN}{\vec \xi}
\newcommand{\nN}{\vec n}
\newcommand{\xN}{\vec x}
\newcommand{\thetaN}{\vec \theta}
\renewcommand{\d}[0]{{\mathrm{d}}}
\newcommand{\Rt}{R_{\theta_1}}
\newcommand{\RT}{R^T_{\theta_1}}
\newcommand{\Pout}[0]{P_{out}^{(1)}}
\newcommand{\tP}{{\tilde P}}
\newcommand{\tPo}[1]{\tilde P^{(#1)}}
\newcommand{\tx}{{\tilde \x}}
\newcommand{\txo}{{\tilde \x_1}}
\renewcommand{\tt}{{\tilde \theta}}
\newcommand{\tto}{{\tilde \theta_1}}
\newcommand{\txi}{\tilde \xi}
\newcommand{\txio}{\tilde \xi_1}
\newcommand{\rhombus}{\mathcal R_\tt}
\newcommand{\e}{{\bf e}}
\newcommand{\im}{\imath}

\usepackage{amsopn}
\DeclareMathOperator{\diag}{diag}
\DeclareMathOperator*{\Res}{Res}

\newcommand{\R}[0]{\mathbb{R}}
\newcommand{\N}[0]{\mathbb{N}}

\newtheorem{remark}{\bf Remark}[section]

\begin{document}

\title[Derivation of a macroscopic model for Brownian hard needles]{Derivation of a macroscopic model for Brownian hard needles}

\author{M. Bruna$^{1}$}
\author{S. J. Chapman$^{2}$}
\author{M. Schmidtchen$^{3}$}

\address{$^{1}$Department of Applied Mathematics and Theoretical Physics, University of Cambridge, Wilberforce Road, Cambridge CB3 0WA, United Kingdom}
\address{$^{2}$Mathematical Institute, University of Oxford, Andrew Wiles Building, Radcliffe Observatory Quarter, Woodstock Road, Oxford OX2 6GG, United Kingdom}
\address{$^{3}$Institute of Scientific Computing, Technische Universit\"at Dresden, Zellescher Weg 25, 01217 Dresden, Germany}
\maketitle

\begin{abstract}
We study the role of anisotropic steric interactions in a system of hard Brownian needles. Despite having no volume, non-overlapping needles exclude a volume in configuration space that influences the macroscopic evolution of the system. 
Starting from the stochastic particle system, we use the method of matched asymptotic expansions and conformal mapping to systematically derive a nonlinear nonlocal partial differential equation for the evolution of the population density in position and orientation.  
We consider the regime of high rotational diffusion, resulting in an equation for the spatial density that allows us to compare the effective excluded volume of a hard-needles system with that of a hard-spheres system. We further consider spatially homogeneous solutions and find an isotropic to nematic transition as density increases, consistent with Onsager's theory.
\end{abstract}
\vskip .4cm
\begin{flushleft}
    \noindent{\makebox[1in]\hrulefill}
\end{flushleft}
	2020 \textit{Mathematics Subject Classification.} 35C20, 35K55, 35Q84,  60J70, 82C22, 70K20, 35B36
	\newline\textit{Keywords and phrases.} many-particle systems, anisotropic particles, excluded-volume interactions, phase transitions, coarse-graining\\[-2.em]
\begin{flushright}
    \noindent{\makebox[1in]\hrulefill}
\end{flushright}
\vskip 1.5cm

\section{Introduction} \label{sec:intro}

Systems of interacting particles are ubiquitous in nature. Examples include biomolecules (e.g. proteins), polymers (e.g. DNA), cells (e.g. bacteria), all the way to multi-cellular organisms (animals). Interactions between organisms may be attractive (keeping a herd cohesive), aligning (keeping animals moving in the same direction), or repulsive (keeping particles a safe distance apart) \cite{DOrsogna:2006ci}. Short-ranged repulsive interactions with singular or hard-core potentials are used to model steric or excluded-volume interactions \cite{BC12}.

Anisotropy plays a crucial role in self-organisation. For example, the helical form of the DNA strand is due to highly anisotropic interactions between DNA bases \cite{snodin2015introducing}. The molecular shape of liquid crystals leads to their remarkable properties \cite{bahadur1990liquid}. Self-propulsion in active matter systems can lead to motility-induced phase separation \cite{Cates:2014tr}, where the uniform suspension becomes unstable and dense clusters of almost stationary particles emerge \cite{Cates:2014tr}. Alignment interactions have been shown to explain the emergence of flocking and milling \cite{carrillo2010particle}. 

Tools to study the rich collective properties of such systems range from simulations at the microscopic level (e.g. molecular dynamics or Monte Carlo simulations) to the study of macroscopic models for statistical quantities, often involving partial differential equations (PDEs). While microscopic models provide a detailed system description, simulating them can become computationally prohibitive. This is due to the large number of particles and the complexity of interactions often involved, mainly if one is after statistical properties (which require averaging over multiple simulations). Macroscopic models operate at the statistical level and can often provide the insight lacking in their microscopic counterparts.

Anisotropy in particle systems comes in many forms. Models can be classified into either first- or second-order models and either soft- or hard-core interactions. 
 In second-order models (which track particles' positions and velocities), particles may interact differently depending on their relative velocities. 
Examples with weak interactions include the Cucker--Smale model \cite{Cucker:2007gt}, and the Vicsek model \cite{Vicsek.1995}, which include alignment interactions in velocities.
One may also add a cone of vision such that an individual only aligns velocity with neighbours within the cone \cite{carrillo2010particle}.
The Cucker--Smale and Vicsek models, and their many variants, have been the starting point in multiple works concerned with deriving kinetic PDE models starting from such microscopic dynamics. It is common to consider a weak or mean-field scaling $1/N$ of the interactions (where $N$ is the number of particles), leading to nonlocal and nonlinear kinetic PDE in the limit \cite{Degond:2017ha,jabin2017mean, topaz2006nonlocal,hoffmann2017keller,carrillo2009double}. The focus in most kinetic models is on how the interaction rule depends on relative positions and velocities, not particles' shapes. An exception is the recent works \cite{HKMS2020, KS2021}, where they consider a system of kinetic hard needles that align upon collision. Instead of a mean-field scaling, they consider the Boltzmann--Grad limit of infrequent but strong interactions and invoke propagation of chaos to derive a closed kinetic equation.

First-order models for anisotropic particles often consider the particle position and orientation and assume diffusive behaviour in both. 
For isotropic particles, microscopic models are well-established (the hard-sphere, the Lennard--Jones, the Coulomb potential, etc.), and current efforts primarily focus on deriving macroscopic models from them. In contrast, anisotropic particles are much harder to model, even at the microscopic scale. There is a trade-off between the complexity of the particles' shape, on the one hand, and the model's analytic tractability, on the other. Interactions may be soft or hard; depending on the application, either may be seen as the `true dynamics'. For example, while soft interactions may be more appropriate for molecules, hard steric interactions may be more fitting for cells, bacteria and animals. 

The most well-known soft anisotropic potential is the Gay--Berne potential \cite{GB1981}. It builds on the work of Berne and Pechukas \cite{BP1972}, which proposed to represent particles as a union of Gaussian potentials and their interaction as the overlap integral of their Gaussians. The Gay--Berne potential combines this anisotropic overlap model with the Lennard--Jones potential (a 12-6 attractive-repulsive potential). 
The multi-phase-field approach \cite{MPV2015, WV2021} is at the other end of the complexity-tractability trade-off. Here, each particle is not characterised by its centre of mass and orientation but by a phase field variable, $\phi_i(x,t) \in [0,1]$, such that $\phi_i(x,t) \approx 1$ if location $x$ is occupied by the particle $i$ at time $t$ and, conversely, $\phi(x,t) \approx 0$ if the $i^{\mathrm{th}}$ particle does not occupy location $x$ at time $t$. 
Due to the diffusive interface between the two states (occupied and unoccupied), repulsive interactions are incorporated in a fashion similar to that of Berne and Pechukas: the overlap integral between two particles (now represented by phase-field variables rather than Gaussians) is computed, and the evolution is such that it minimises the area of overlap. 

\begin{figure}[ht]
\begin{center}
\input{excluded_region.pstex_t} 
\includegraphics[width=0.5\textwidth]{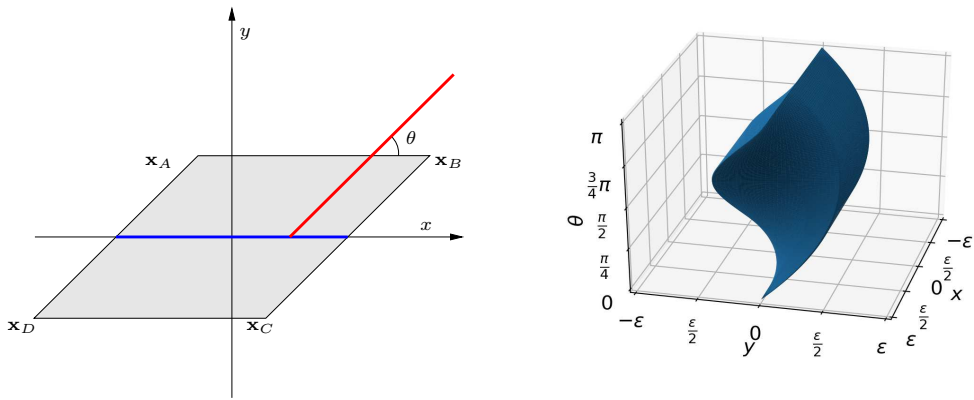}
\end{center}
\caption{Left panel: excluded volume of a horizontal hard needle (blue) with centre at the origin. The centre of a second needle (red) with orientation $\theta$ cannot be placed inside the excluded region (shaded grey area) as it would lead to an overlap. Right panel: excluded volume in phase space. The vertical axis denotes the relative angle between two needles.}
\label{fig:excluded_volume}
\end{figure}

The models above have in common that the space taken up by particles is not precisely localised, in contrast to hard-core models. 
Hard-core ellipsoids and rods are the natural generalisations of hard spheres to model anisotropy. A hard-core particle induces an excluded region (where no other particle can enter). In his seminal paper \cite{Onsager:1949jk}, Onsager 
finds expressions for the excluded volume of various particle shapes such as ellipses, discs, and rods. 
The most striking example in his treatise is a hard needle of length $\epsilon$ in two dimensions, which has zero volume but excludes a volume in configuration space of $\epsilon^2 |\sin(\theta)|$ to another needle with relative orientation $\theta$ (see Figure \ref{fig:excluded_volume}). The problem of interacting needles in three dimensions is fundamentally different since needles exclude zero volume in configuration space in addition to having no volume.
Rods with a core-shell structure have also been used in microscopic models of self-assembly \cite{Farrell:2012ks} and morphogenesis in bacterial colonies \cite{Doumic:2020uo}. 
Interesting mathematical problems arise from considering even just one anisotropic hard-core particle. For example, in  \cite{Holcman:2012dz}, they study the mean turnaround time of a Brownian needle in a narrow planar strip as a simplified model for mRNA or stiff DNA fragments under extreme confinement. In \cite{chen2021shape}, they consider an anisotropic Brownian microswimmer in a channel and show that no-flux boundary conditions with the flat channel walls lead to nontrivial boundaries in configuration space.  

Since one hard-core anisotropic particle already poses mathematical challenges and, ordinarily, natural systems comprise large ensembles of anisotropic particles, it is easy to see that their study is substantially more challenging.
This explains the dearth of macroscopic PDE models systematically derived from underlying dynamics and the popularity of computational and phenomenological approaches to incorporate anisotropic interactions in PDE models.
Phenomenological models have been widely used in the context of polymer and liquid crystals theory. These include the so-called tube theory \cite{doi1986theory}, which assumes polymers as rigid filaments that, under crowding, move along a tube formed by the surrounding polymers, as well as the Landau-de Gennes $Q$-tensor theory for nematic liquid crystals \cite{de1993physics}, which represents polar molecules via a continuum order parameter. A lot of work has been dedicated to validating these theories by comparing their predictions with microscopic models with different levels of success \cite{FM1981}. 

In \cite{MKFL2020}, they consider a system of self-propelled needles, with collisions such that energy and momentum are preserved, and validate the tube theory (they find that the self-diffusion coefficient of a needle increases with concentration, in contrast to that of hard spheres \cite{Bruna:2012wu}).

In this paper, we focus on, possibly, the simplest hard-core anisotropic system, namely that of $N$ Brownian needles of length $\epsilon$ with non-overlapping constraints in two dimensions. Using matched asymptotic expansions, we systematically derive a macroscopic PDE model in the asymptotic regime $\epsilon^2 N \ll 1$. To our knowledge, this is the first systematic derivation for such a system. We take an approach similar to \cite{BC12}, in which the authors consider a system of $N$ Brownian hard disks of diameter $\epsilon$ under a drift $\f(\x)$ in two spatial dimensions. Under the assumption that the volume fraction of the particles is small,  the one-particle probability density $\rho(\x,t)$ satisfies the nonlinear diffusion equation (Eq. (11) of Ref. \cite{BC12})
\begin{align}
\label{al:brunachapman_equation}
\frac{\partial \rho}{\partial t}(\x, t) = 
	\nabla_\x \cdot \left \{ [1+ \pi(N-1)\epsilon^2 \rho] \nabla_\x \rho - \f(\x) \rho \right \},
\end{align}
in $\R^2$. The goal of this paper is to derive an analogous PDE to \eqref{al:brunachapman_equation} for the one-particle density $p(\x,\theta,t)$ describing the probability of a needle with centre at $\x$ and orientation $\theta$ at time $t$.

The structure of the paper is as follows. In Section \ref{sec:the_model}, we introduce the particle-based model, that is, a system of $N$ Brownian needles with drifts and its associated Fokker-Planck equation describing the whole ensemble probabilistically. Section \ref{sec:derivation_macro} is devoted to the systematic derivation of the effective model using the method of matched asymptotic expansions. Section \ref{sec:high_rotation} is dedicated to systems with high rotational diffusion coefficients, which are shown to inherit striking similarities with the hard-disk model, Section \ref{sec:space_hom}, proposed by \cite{BC12}. We conclude the paper with the space-homogeneous model and, upon performing a linear stability analysis, we find that the system exhibits an isotropic-to-nematic phase transition consistent with Onsager's theory \cite{PhysRevA.17.2067}.

\section{The microscopic model and its associated Fokker--Planck equation} \label{sec:the_model}

We start by describing the individual-based (microscopic) model. We suppose there are $N\in\N$ identical hard needles of length $\epsilon$ distributed in a bounded domain $\Omega\in\R^2$. For $1\leq i \leq N$, we denote by $\X_i(t)\in \Omega$ the centre of the $i$th needle and by $\Theta_i(t) \in [0,\pi]$  its orientation. We choose $\Omega$ to be the 2-dimensional torus $\mathbb T = \R/(\pi \mathbb Z)$, imposing periodic boundary conditions. The $\pi$-period is chosen for mathematical convenience, such that $\Upsilon = \Omega \times [0,\pi) = \mathbb T^3$.
 As pointed out in the introduction, the spatial extension of the needles restricts their ability to move freely in the domain due to non-overlapping constraints, in contrast to a system of point particles.

Each needle evolves according to a translational (resp. rotational) Brownian motion with diffusion constant $D_T$ (resp. $D_R$) in an external force field $f = (\f_T, f_R)$ that may depend both on the position and orientation of the needle, but not on the other needles.
This leads to the system of stochastic differential equations (SDEs)
\begin{subequations} \label{sdes}
	\begin{align} \label{sde_x}
	\d \X_i(t) &= \sqrt{2 D_T} \, \d {\bf W}_{T,i}(t) + \f_T(\X_i, \Theta_i) \d t,\\
	\label{sde_theta}
	\d \Theta_i(t) &= \sqrt{2 D_R} \, \d W_{R,i}(t) + f_R(\X_i, \Theta_i)\d t,
\end{align}
\end{subequations}
for $1 \le i \le N$ and $(\X_i, \Theta_i) \in \Upsilon$. Here ${\bf W}_{T,i}$ and $W_{R,i}$ are standard independent Brownian motions for $1 \le i \le N$. 
In addition, we impose reflective boundary conditions whenever two needles come into contact, thereby introducing a coupling to an otherwise uncoupled system of $N$ SDEs. 

It is convenient to consider the  \emph{joint probability density} $P_N(\xiN, t)$ associated to system \eqref{sdes}, where $\xiN = (\xi_1, \dots, \xi_N)$ and $\xi_i = (\x_i, \theta_i)$, for $1\le i\le N$. The density $P_N$ describes the probability of the entire system of $N$ needles being in state $\xiN$ at time $t$. It is well known that $P_N$ satisfies the Fokker--Planck equation
\begin{subequations}
\label{eq:highdim-fokker-planck}
	\begin{equation}
		\label{al:fokker_planck_for_needle_system}
		\partial_t P_N =  \nabla_{\xN} \cdot \left [ D_T \nabla_{\xN} P_N - \vec F_T(\xiN) P_N \right ]  + \nabla_{\thetaN} \cdot \left [D_R \nabla_{\thetaN} P_N - \vec F_R(\xiN) P_N \right ],
	\end{equation}
where $\xN = (\x_1, \dots, \x_N)$, $\thetaN = (\theta_1, \dots, \theta_N)$, $\vec F_T(\xiN) = (\f_T(\xi_1), \dots, \f_T(\xi_N))$ and $\vec F_R(\xiN) = (f_R(\xi_1), \dots,\allowbreak f_R (\xi_N))$. Due to the hard-core interactions between needles, note that \eqref{al:fokker_planck_for_needle_system} is not defined on $\xiN \in \Upsilon^N$ but its perforated form $\Upsilon_{\epsilon}^N :=\Upsilon^N \setminus \mathcal{B}_\epsilon^N$. Here, $\mathcal{B}_{\epsilon}^N$ denotes the set of  \emph{illegal configurations} where at least two needles overlap, i.e., 
\begin{align*}
	\mathcal{B}_{\epsilon}^N:=\big\{\xiN \in \left(\Upsilon \right)^N \ | \ \exists i\neq j \mbox{ s.t. }  \mathcal{N}(\xi_i)\cap\mathcal{N}(\xi_j) \neq \emptyset \big\},
\end{align*}
where
$$
\mathcal{N}(\x,\theta):= \left\{
	x + \lambda
	\begin{pmatrix}
			\cos(\theta)\\
			\sin(\theta)
	\end{pmatrix}
	\bigg|  \ |\lambda| \leq \frac\epsilon2
	\right\},
$$
denotes the set of all points belonging to a needle at $(\x,\theta)$. 

On $\partial \Upsilon_{\epsilon}^N$ (corresponding to configurations with at least two needles in contact), we prescribe reflective boundary conditions
\begin{equation}
\label{al:fokker_planck_for_needle_system_boundary_condition}
\begin{bmatrix} 
D_T \nabla_{\xN} P_N - \vec F_T(\xiN) P_N \\
D_R \nabla_{\thetaN} P_N - F_R(\xiN) P_N
\end{bmatrix}\cdot \nN = 0, \quad \text{on} \quad \partial \Upsilon_{\epsilon}^N,
\end{equation}
\end{subequations}
where $\nN \in \mathcal{S}^{3N-1}$ denotes the unit outward normal on the boundary. Finally, we assume that the initial positions of the particles are identically distributed so that the initial condition $P(\xiN, 0) = P_0(\xiN)$ is invariant to permutations of the particles' labels. 

\section{Derivation of the macroscopic model} \label{sec:derivation_macro}
In the previous section, we have established a connection between the particle-based dynamical system \eqref{sdes} and the associated Fokker-Planck equation \eqref{eq:highdim-fokker-planck}. We highlight that the dimensionality of both descriptions increases as more needles are added to the system, rendering their analytical or numerical study intractable. 
This section is dedicated to deriving an effective model in the form of a nonlinear evolution equation for the \emph{one-particle probability density}
\begin{equation}
\label{eq:firstmarginal}
	p(\xi, t) := \int_{\Upsilon_\epsilon^N} P_N(\vec \xi,t) \delta(\xi_1 - \xi) d \vec \xi.
\end{equation}

In the case of $\epsilon = 0$, the needles become point particles and, as a consequence, their evolutions \eqref{sdes} decouple and, for suitable iid initial conditions, we have that  
$$
	P_N(\xiN, t) =  \prod_{i=1}^N p(\xi_i,t).
$$
In this setting, the first marginal is shown to satisfy the following equation
\begin{equation}
	\label{points}
		\partial_t p(\xi, t) =  \nabla_{\x} \cdot \left [ D_T \nabla_{\x} p - \f_T(\xi) p \right ]  + \partial_{\theta} \left [D_R \partial_{\theta} p - f_R(\xi) p \right ],
\end{equation}
with $t\ge 0$ and $\xi \in \Upsilon$. 
Unlike point particles, needles of length $\epsilon>0$ exclude a certain volume in phase space. 
\begin{remark}[Excluded region of a needle] \label{rem:excluded_region}
The region in phase space excluded by a needle at $\xi_1$ is denoted by 	$B_\epsilon(\xi_1)$ (Fig.~\ref{fig:excluded_volume}). 
Depending on the relative orientation $\theta := \theta_2-\theta_1$ between the two needles, the cross-section of $B_\epsilon$ for fixed $\theta$ range from a line of length $2 \epsilon$  ($\theta = 0$) to a square of side $\epsilon$ ($\theta = \pi/2$). For general $\theta$, the slice is a rhombus of area $\epsilon^2 \sin \theta$ with nodes at
\begin{equation}
	\label{vertices}
\begin{alignedat}{2}
\x_A &= \x_1 + \frac{\epsilon}{2} \Rt (-1 + \cos \theta, \sin \theta),	 &\quad \x_B &= \x_1 + \frac{\epsilon}{2} \Rt(1 + \cos \theta, \sin \theta), \\
\x_C &= \x_1 + \frac{\epsilon}{2} \Rt(1 - \cos \theta, -\sin \theta), &\quad \x_D &= \x_1 + \frac{\epsilon}{2} \Rt(-1 - \cos \theta, -\sin \theta),
\end{alignedat}
\end{equation}
where $\Rt$ is the rotation matrix
\begin{equation} \label{rotation_matrix}
	\Rt = \begin{pmatrix}
		\cos \theta_1 & -\sin \theta_1  \\
		\sin \theta_1 & \cos \theta_1
	\end{pmatrix}.
\end{equation}
We denote by $\hat n_2$ the outward unit normal on $B_\epsilon(\xi_1)$ (outward of $\Upsilon(\xi_1)$ so it points into the shaded area in Fig.~\ref{fig:excluded_volume}). If the boundary of $B_\epsilon(\xi_1)$ is given by the relation $\chi(\xi_2) = 0$, we have that $\hat n_2 \propto \nabla_{\xi_2} \chi $. 
For example, the top edge $\x_A \x_B$ is given by $\chi(\xi_2) = y_A + \tan \theta_1 (x_2 - x_A) - y_2$ = 0 and the normal vector is 
\begin{equation} \label{normal_xDxA}
	\hat n_2 \propto \nabla_{\xi_2} \chi = \left(\tan \theta_1, 1, 
\frac{\epsilon}{2} ( \cos \theta_2 + \tan \theta_1 \sin \theta_2) \right).
\end{equation}
\end{remark}

For $\epsilon>0$, the equation for the one-particle density $p(\xi_1, t)$ is obtained by integrating \eqref{eq:highdim-fokker-planck} with respect to $\xi_2, \dots, \xi_N$ for $\xi_1$ fixed. The perforations in $\Upsilon_\epsilon^N$ lead to boundary integrals for $\xi_i\in B_\epsilon(\xi_1)$ on which the two-particle probability density $P_2(\xi_1, \xi_i, t)$ needs to be evaluated. One can go back to \eqref{eq:highdim-fokker-planck} and obtain an equation for $P_2$, which in turn depends on the three-particle probability density $P_3$. This is known as the BBGKY hierarchy. In this work, we assume that $ \phi = \epsilon ^2 N \ll 1$ such that this hierarchy can be truncated ``asymptotically''. 

We note from Remark \ref{rem:excluded_region} that the volume of $B_\epsilon(\xi_1)$ is $\epsilon^2 \int_0^\pi \sin(\theta) \d \theta = 2 \epsilon^2$. If $\phi \ll 1$, the volume in $\Upsilon_\epsilon^N$ occupied by configurations where two needles are closeby is $O(\phi)$, whereas the volume of configurations where three or more needles are nearby is much smaller ($O(\phi^2)$). 
Hence, it means that, at the leading order, the equation for $p$ coincides with the point particles equation \eqref{points} and that the first correction appears at $O(\phi)$ and is due to two-needle interactions. Three- and more-needle interactions are higher-order corrections. 
Therefore, we may neglect three-particle interactions in the equation for $P_2(\xi_1,\xi_2, t)$ and consider 
\begin{subequations}
	\label{fp_2}
\begin{align}
\label{eqn:fokker_planck_two_needles}
	\partial_t P_2 = \nabla_{\xi_1} \cdot \left[D \nabla_{\xi_1} P_2 -f(\xi_1) P_2 \right] + \nabla_{\xi_2} \cdot \left[D \nabla_{\xi_2} P_2 -f(\xi_2) P_2 \right],
\end{align}
in $\Upsilon_\epsilon^2$, where $D = \diag (D_T, D_T, D_R)$ and $f(\xi) = (\f_T(\xi), f_R(\xi))$, together with reflecting boundary conditions
\begin{equation}
\label{eqn:fokker_planck_two_needles_boundary_condition}
\left[D \nabla_{\xi_1} P_2 - f(\xi_1) P_2 \right] \cdot n_1 + \left[
D \nabla_{\xi_2} P_2 - f(\xi_2) P_2 \right] \cdot n_2 = 0 
\end{equation}
on $\partial \Upsilon_\epsilon^2$. Here $n_1$ (resp. $n_2$) are the components of the unit normal $\vec n$ corresponding to the coordinates of the first (resp. second) needle. It turns out that $n_1 = -n_2$ such that
$\vec n =  \sqrt 2/2 (- \hat n_2,  \hat n_2)$, where $\hat n_2$ is defined in Remark \ref{rem:excluded_region}. 
\end{subequations}

\subsection{Evolution of the first marginal}

Let $\Upsilon(\xi_1) = \Upsilon \setminus B_\epsilon(\xi_1)$ denote the second particle's phase space given that the first particle is in state $\xi_1$. Integrating \eqref{eqn:fokker_planck_two_needles} over $\Upsilon(\xi_1)$ yields
\begin{align}
\begin{split}
	\label{eq:intermediate_evolution_equation}
\partial_t p(\xi_1, t) &= \int_{\Upsilon(\xi_1)} \partial_t P_2(\xi_1, \xi_2, t) \d \xi_2 \\
	&= \int_{\Upsilon(\xi_1)} \nabla_{\xi_1}\cdot \left[ D \nabla_{\xi_1} P_2 - f(\xi_1) P_2 \right] \d \xi_2 + \int_{\Upsilon(\xi_1)} \nabla_{\xi_2} \cdot \left[ D \nabla_{\xi_2}P_2 -f(\xi_2) P_2 \right]\d \xi_2.
	\end{split}
\end{align}
Using Reynold's transport theorem, the first integral becomes
\begin{align}
\begin{split}
	\label{eq:equation_for_A}
\int_{\Upsilon(\xi_1)} &\nabla_{\xi_1}\cdot \left[ D \nabla_{\xi_1} P_2 - f(\xi_1) P_2 \right] \d \xi_2 \\
&=\nabla_{\xi_1} \cdot \left[ D \nabla_{\xi_1} p-f (\xi_1) p \right] 
	+ \oint_{\partial B_\epsilon(\xi_1)} \left[ f(\xi_1) P_2 - 2 D \nabla_{\xi_1}P_2 - D \nabla_{\xi_2} P_2 \right] \cdot \hat n_2 \d S_{\xi_2},
\end{split}
\end{align}
The second integral in \eqref{eq:intermediate_evolution_equation} is
\begin{align}
	\label{eq:equation_for_B}
\int_{\Upsilon(\xi_1)} \nabla_{\xi_2} \cdot \left[ D \nabla_{\xi_2}P_2 -f(\xi_2) P_2 \right]\d \xi_2 = \oint_{\partial B_\epsilon(\xi_1)} \left[ D \nabla_{\xi_2} P_2  - f (\xi_2)P_2 \right] \cdot \hat n_2 \d S_{\xi_2}.
\end{align}
Substituting \eqref{eq:equation_for_A} and \eqref{eq:equation_for_B} into  \eqref{eq:intermediate_evolution_equation} we obtain
\begin{align}
	\label{eq:evolution_eqn_with_bdry_integral}
	\partial_t p(\xi_1, t) = \ &\nabla_{\xi_1} \cdot \left[ D \nabla_{\xi_1} p - f(\xi_1) p \right] + I,
\end{align}
where the collision integral $I$ is
\begin{equation}
	\label{col_integral2}
	I = - \oint_{\partial B_\epsilon(\xi_1)}   D  \left( \nabla_{\xi_1} P_2 + \nabla_{\xi_2} P_2 \right) \cdot \hat n_2 \d S_{\xi_2}.
\end{equation}

The evolution equation \eqref{eq:evolution_eqn_with_bdry_integral} for the first marginal $p$ still depends on the joint probability density function $P_2$. A common approach to overcome this is to use a closure assumption, for instance, the mean-field approximation, $P_2(\xi_1,\xi_2,t)=p(\xi_1,t)p(\xi_2,t)$. 
However, such an approach ignores correlations between both particles, and it is not suitable for systems of strongly interacting particles with short-range repulsive interactions such as hard needles. Instead, we employ the method of matched asymptotics to compute the collision integral $I$ systematically. 

\subsection{Matched asymptotics expansions}

We introduce a partition of the domain $\Upsilon(\xi_1)$ consisting of an \emph{inner region}, when the two needles are close to each other, $\|x_1 - x_2\|_2 \sim \epsilon$, and an \emph{outer region}, when the two needles are far apart, ${\|\x_1 - \x_2\| \gg \epsilon}$. In the outer region, we suppose that particles are independent at leading order, whereas we consider their correlation in the inner region. 

In the outer region we define $P_{out}(\xi_1, \xi_2, t) = P_2(\xi_1, \xi_2, t)$. Then by independence, the two-particle density function is\footnote{Independence only tells us that $P_{out}(\xi_1, \xi_2, t) \sim q(\xi_1, t) q(\xi_2, t)$ for some function $q$, but the normalisation condition on $P_2$ implies $p = q + O(\epsilon)$.}
\begin{equation}
	\label{outer_ansatz}
	P_{out}(\xi_1, \xi_2, t) = p(\xi_1, t) p(\xi_2, t) + \epsilon \Pout (\xi_1, \xi_2, t) + \cdots.
\end{equation}

In the inner region, we introduce the \emph{inner variables} $\txio = (\txo, \tto)$ and $\txi = (\tx, \tt)$, defined as
\begin{equation}
	\label{change_vars}
\begin{alignedat}{2}
	\x_1 &= \txo, & \qquad \x_2 & = \txo + \epsilon \Rt \tx,\\
	\theta_1 &= \tto, &  \theta_2 & = \tto +  \tt,
\end{alignedat}
\end{equation}
and the inner function $\tP (\txio, \txi,t) = P_2(\xi_1, \xi_2, t)$. 
The coordinates $(\tx, \tt)$ define the configuration of the second needle relative to the first. The excluded volume $B_\epsilon(\xi_1)$ becomes $B_1(0)$ in inner variables. In the $\txi$-space, this is now a volume centred at the origin with two horizontal sides (see Figure \ref{fig:excluded_volume} and Remark \ref{rem:excluded_region}). 
Using that $\tx = \epsilon^{-1} \RT (\x_2 - \x_1)$, the  derivatives transform according to
\begin{equation*}
	\begin{alignedat}{2}
		\nabla_{\x_1} &\to \nabla_\txo - \epsilon^{-1} \Rt\nabla_\tx, & \qquad \nabla_{\x_2} &\to \epsilon^{-1} \Rt \nabla_\tx,\\
		\partial_{\theta_1} &\to \partial_{\tto} - \partial_\tt + \tilde y \partial_{\tilde x} - \tilde x \partial_{\tilde y}, & \partial_{\theta_2} &\to \partial_\tt.
	\end{alignedat}
\end{equation*}

In terms of the inner variables, \eqref{eqn:fokker_planck_two_needles} reads
\begin{align}
\label{fp_inner}
\begin{split}
	\epsilon^2 \partial_t\tP  
	= \ & 2D_T \Delta_\tx \tP \\
	&- \epsilon
	    \left(
    		2D_T \nabla_\txo\cdot \Rt\nabla_\tx \tP  +  \nabla_\tx \cdot \left\{ \Rt\big[ \f_T(\txo + \epsilon \tx, \tto + \tt) - \f_T(\txio) \big]\tP \right\}
		\right)\\
		&+\epsilon^2 \bigg\{ \nabla_\txo \cdot 
			\left[ D_T \nabla_\txo \tP - \f_T(\txio) \tP \right] + D_R \left[ \left(
					\partial_\tto  - \partial_{\tt} + \tilde y \partial_{\tilde x} - \tilde x \partial_{\tilde y}\right)^2 + \partial_\tt^2 \right] \tP\\
				&\qquad \quad - (\partial_{\tto} - \partial_\tt + \tilde y \partial_{\tilde x} - \tilde x \partial_{\tilde y}) \left[ f_R(\txio) \tP \right] - \partial_\tt \left[ f_R(\txo + \epsilon\tx, \tto + \tt) \tP \right]
		\bigg\}.
\end{split}
\end{align}
In order to write the boundary condition \eqref{eqn:fokker_planck_two_needles_boundary_condition} 
in terms of the inner variables, we need to determine how the normal $\hat n_2$ changes under the transformation. Following the procedure in Remark \ref{rem:excluded_region}, we have $\nabla_{\xi_2} \chi \to (\epsilon^{-1} \Rt \nabla_{\tx} \tilde \chi, \partial_\tt \tilde \chi)$, where $\tilde \chi(\tilde \xi) = 0$ describes the boundary in inner variables. Therefore
\begin{equation} \label{rem:normal_transform}
	\hat n_2 \to (\Rt \tilde \n, \epsilon \tilde n_\theta).
\end{equation}
For example, the top edge $\x_A \x_B$ becomes $\tilde \chi =  \sin(\tt) - \tilde y = 0$ and the normal vector in the inner variables is
$
		\tilde n = (\tilde \n, \tilde n_\theta) \propto \nabla_{\tx} \tilde \chi = (0, -1, \cos \tt).
$
Using \eqref{rem:normal_transform} and $n_1 = -n_2$ as pointed out earlier, the no-flux boundary condition \eqref{eqn:fokker_planck_two_needles_boundary_condition} becomes 
\begin{align} \label{bc_inner}
\begin{split}
		0 &= \left \{ 2 D_T \Rt \nabla_{\tx} \tP - \epsilon D_T \nabla_{\txo} P - \epsilon  \left[  \f_T(\txo + \epsilon\tx, \tto + \tt) - \f_T(\txio) \right]  \tP \right \} \cdot \Rt \tilde \n \\
		& + \epsilon^2 \left \{  D_R \big[ 2\partial_{\tilde \theta} \tP -  \partial_{\tilde \theta_1} \tP + \tilde x \partial_{\tilde y} \tP - \tilde y \partial_{\tilde x} \tP \big] - \big[f_R(\txo + \epsilon\tx, \tto + \tt)  - f_R(\txio) \big] \tP \right \}  \tilde n_\theta,
\end{split}
\end{align}
for $\txi \in \partial B_1(0)$. 
Finally, we impose the \emph{matching boundary condition} to ensure that, as the two needles become further apart and enter the outer region, the inner solution $\tP$ will match with the outer solution $P_{out}$. Expanding \eqref{outer_ansatz} in the inner variables, 
\begin{align} \label{matching}
\begin{split}
	\tP &\sim p(\txo,\tto,t) p(\txo + \epsilon \Rt \tx, \tto + \tt,t) + \epsilon \Pout (\txo,\tto, \txo + \epsilon \Rt \tx, \tto + \tt,t) \\
	& \sim p p ^+ + \epsilon \left[ p \Rt\tx \cdot \nabla_\txo p^+ + \Pout (\txo,\tto, \txo, \tto + \tt,t) \right], \qquad \text{as} \qquad |\tx| \to \infty, 
\end{split}	
\end{align}
where $p := p(\txo,\tto, t)$ and $p^+ := p(\txo, \tto + \tt, t)$. 

We look for a solution of \eqref{fp_inner}, \eqref{bc_inner}, and \eqref{matching} of the form  $\tP = \tPo{0} +  \epsilon \tPo{1} + \cdots$. The leading-order problem is 
\begin{equation} 
\label{inner_0}
    \left\{
        \begin{alignedat}{2}
            \Delta_\tx \tPo{0}	 & = 0,  & &\\
            \Rt\nabla_\tx \tPo{0} \cdot \Rt\tilde \n & = 0,   &\qquad &\txi \in \partial B_1(0),\\
            \tPo{0} & \sim p p^+, & &  |\tx| \sim \infty.
        \end{alignedat}
    \right.
\end{equation}
This is a problem in the inner spatial variables $\tx$, and that $\txo$, $\tto$, and $\tt$ can be regarded as parameters. In particular, \eqref{inner_0} is defined for $\tx \in \mathbb R^2 \setminus \rhombus$, where $\rhombus$ denotes the rhombus corresponding to slicing the excluded volume $B_1(0)$ at $\tt$ (see Figure \ref{fig:excluded_volume}). The solution of \eqref{inner_0} is
\begin{equation}
	\label{P0_sol}
	\tPo{0} = p p^+.
\end{equation}
 Using \eqref{P0_sol} and expanding $\f_T$, the $O(\epsilon)$ problem reads
\begin{equation} \label{inner_1}
	\begin{alignedat}{2}
\Delta_\tx \tPo{1}	 & = 0,  & & \tx \in \mathbb R^2 \setminus \rhombus, \\
\Rt \nabla_{\tx} \tPo{1} \cdot \Rt \tilde \n &=  \frac{1}{2} \left[  \nabla_{\txo} (p p^+) + \frac{(p p^+)}{D_T} (\f_T^+ - \f_T)   \right ] \cdot \Rt \tilde \n,   &\qquad &\tx \in \partial \rhombus,\\
\tPo{1} & \sim p \nabla_\txo p^+ \cdot \Rt \tx  + \Pout (\txo,\tto, \txo, \tto + \tt, t), & &  |\tx| \sim \infty,
\end{alignedat}
\end{equation}
where $\f_T := \f_T(\txo,\tto)$ and $\f_T^+ := \f_T(\txo, \tto + \tt)$. 
We can rewrite problem \eqref{inner_1} as 
\begin{equation} 
    \label{inner_1b}
    \left\{
    \begin{alignedat}{2}
        \Delta_\tx \tPo{1}	 & = 0,  & & \tx \in \mathbb R^2 \setminus \rhombus, \\
        \nabla_{\tx} \tPo{1} \cdot \tilde \n &=  \RT{\bf A} \cdot \tilde \n,   &\qquad &\tx \in \partial \rhombus,\\
        \tPo{1} & \sim \RT{\bf B}_\infty \cdot \tx + C_\infty, & &  |\tx| \sim \infty,
    \end{alignedat}
    \right.
\end{equation}
where $\bf A$, ${\bf B}_\infty$ and $C_\infty$ are functions of $\txo, \tto, \tt$ and $t$ only and given by 
\begin{align} \label{ABC_constants}
\begin{aligned}
	{\bf A} &= \frac{1}{2}  \left[  \nabla_{\txo} (p p ^+) + \frac{p p^+}{D_T} (\f_T^+ - \f_T)   \right ], \\ 
	{\bf B}_\infty &= p  \nabla_\txo p^+, \\
	C_\infty &= \Pout (\txo,\tto, \txo, \tto + \tt, t). 
	\end{aligned}
\end{align}
The solution to \eqref{inner_1b} is given by
\begin{equation}\label{p1_decomposed}
\tPo{1} = \RT{\bf A} \cdot \tx + \RT {\bf B} \cdot {\bf u} + C_\infty,	
\end{equation}
where ${\bf B} := {\bf B}_\infty - {\bf A}$ and ${\bf u} = (u_1,u_2)$ satisfy the following problems:
\begin{align} 
    \label{prob_u1}
    \left\{
    \begin{alignedat}{2}
        \Delta_\tx u_1	 & = 0,  & & \tx \in \mathbb R^2 \setminus \rhombus, \\
        \nabla_{\tx} u_1 \cdot \tilde \n &=  0,   &\qquad &\tx \in \partial \rhombus,\\
        u_1 & \sim \tilde x, & &  |\tx| \sim \infty,
    \end{alignedat}
    \right.
\end{align}
and
\begin{equation} 
    \label{prob_u2}
    \left\{
	\begin{alignedat}{2}
        \Delta_\tx u_2	 & = 0,  & & \tx \in \mathbb R^2 \setminus \rhombus, \\
        \nabla_{\tx} u_2 \cdot \tilde \n &=  0,   &\qquad &\tx \in \partial \rhombus,\\
        u_2 & \sim \tilde y, & &  |\tx| \sim \infty.
    \end{alignedat}
    \right.
\end{equation}
Thus we have reduced the inner problem \eqref{inner_1} to two problems for $u_1(\tx)$ and $u_2(\tx)$ that only depend on $\tt$ through their domain of definition,  namely the exterior of a rhombus whose tilting depends on $\tt$ (see Figure \ref{fig:excluded_volume}). Problems \eqref{prob_u1} and \eqref{prob_u2} are solved via conformal mapping in Appendix \ref{sec:conformal}.

\subsection{Collision integral} \label{sec:integral}

In this subsection, we go back to the integrated equation \eqref{eq:evolution_eqn_with_bdry_integral} and use the inner solution $\tP$ to evaluate the collision integral $I$ in \eqref{col_integral2}.
Transforming \eqref{col_integral2} to inner variables, we obtain
\begin{equation}
	\label{col_integral_in}
	I = - \epsilon  D_T \oint_{\partial B_1(0)}  \nabla_{\txo} \tP \cdot \Rt \tilde \n \, \d S_{\txi} - \epsilon^2 D_R \oint_{\partial B_1(0)} \left( \partial_{\tto} \tP + \tilde y \partial_{\tilde x} \tP - \tilde x \partial_{\tilde y} \tP  \right) \tilde n_\tt \, \d S_{\txi},
\end{equation}
using \eqref{change_vars} and \eqref{rem:normal_transform}. 

We evaluate \eqref{col_integral_in} by breaking $I$ in powers of $\epsilon$, $I = I^{(0)} + \epsilon  I^{(1)} + \cdots$. Clearly $ I^{(0)} = 0$. The first-order integral is 
\begin{align*}
	 I^{(1)} = - D_T \oint_{\partial B_1(0)}  \RT \nabla_{\txo} \tPo{0} \cdot \tilde \n \, \d S_{\txi} = 0,
\end{align*}
using that $\tPo{0}$ is independent of $\tx$, see \eqref{P0_sol}, and that we are integrating the normal of a closed curve (for $\tt$ fixed).  
At the next order, we have
\begin{equation} \label{int_o2}
	I^{(2)} = - D_T \underbrace{\oint_{\partial B_1(0)}  \RT \nabla_{\txo} \tPo{1} \cdot \tilde \n \, \d S_{\txi}}_{I_\tx} \ - \ D_R \underbrace{\oint_{\partial B_1(0)} \partial_{\tto} \tPo{0} \tilde n_\tt \, \d S_{\txi} }_{I_\tt},
\end{equation}
using again that $\tPo{0}$ is independent of $\tx$, making the terms $\tilde y \partial_{\tilde x} \tPo{0} - \tilde x \partial_{\tilde y} \tPo{0}$ vanish in the second integral. The latter can be further simplified to
\begin{align}\label{int_theta}
\begin{split}
	I_\tt &=  - \int_{B_1(0)}  \partial_{\tt}  \partial_{ \tto} (p p^+) \, \d \txi = -\partial_{ \tto} \int_{0}^{\pi} \partial_{\tt} (p p^+) \int_{\rhombus} \, \d \tx \d \tt  = - \partial_{ \tto} \int_{0}^{\pi} \sin \tt  \partial_{\tt} (p p^+) \, \d \tt.
	\end{split}
\end{align}
In the first equality, we have applied the divergence theorem to $(0, 0, \partial_{\tto} \tPo{0})$. In the last equality, we have used that $\rhombus$ is the rhombus tilted by angle $\tt$ in inner variables, which has area $\sin \tt$ (see Remark \ref{rem:excluded_region}).
The integral $I_\tx$ in \eqref{int_o2}  can be rewritten as
\begin{equation}
	\label{Ix}
	I_\tx = \int_{0}^{\pi}  \int_{\partial \rhombus}  \RT \nabla_{\txo} \tPo{1} \cdot \tilde \n \, \d S_{\tx} \, \d \tt =  \int_{0}^{\pi}  J(\txo, \tto, \tt) \, \d \tt
\end{equation}
with $J =  \int_{\partial \rhombus}  \RT \nabla_{\txo} \tPo{1} \cdot \tilde \n \, \d S_{\tx}$. Using the expression for $\tPo{1}$ in \eqref{p1_decomposed}, we find that
(see Appendix \ref{sec:integral_all})
\begin{align} 
\label{int_J2}
	\begin{aligned}
	J 
	&= -\nabla_{\txo} \cdot\left ( \sin \tt    {\bf A}  + M(\tto, \tt) {\bf B} \right),
		\end{aligned}
\end{align}
where $M(\tto, \tt) = \Rt T(\tt) \RT$ with $T(\tt)$ the symmetric $2\times 2$ matrix \eqref{matrix_T} whose entries are plotted in Fig.~\ref{fig:Tdiag_Txy}. The matrix $T(\tt)$ is positive definite and contains information on the effect of the excluded volume due to a horizontal needle on a second needle with orientation $\tt$. We observe that:
For $\tt = \pi/2$, the diagonal terms are equal while the cross-terms are zero, as expected, since the excluded region is symmetric (a square). 
For $\tt = 0, \pi$, the needle is ``invisible'' to the horizontal flow ($T_{11} = 0$) and the effect on the vertical flow is maximal ($T_{22}$ largest).
\def \scl {.9}
\begin{figure}[bth]
\unitlength=1cm
\begin{center}
\psfrag{a}[][][\scl]{(a)} \psfrag{b}[][][\scl]{} \psfrag{ai}[][][\scl][-90]{$a_i$} \psfrag{a1}[][][\scl]{$a_1$} \psfrag{a2}[][][\scl]{$a_2$}  
\psfrag{th}[][][\scl]{$\tt$} \psfrag{Txx}[][][\scl]{$T_{11}$} \psfrag{Tyy}[][][\scl]{$T_{22}$}  \psfrag{Txy}[][][\scl]{$T_{12}$}  \psfrag{Tij}[][][\scl][-90]{} \psfrag{Tii}[][][\scl][-90]{$T_{ij}$}
\includegraphics[width=.5\textwidth]{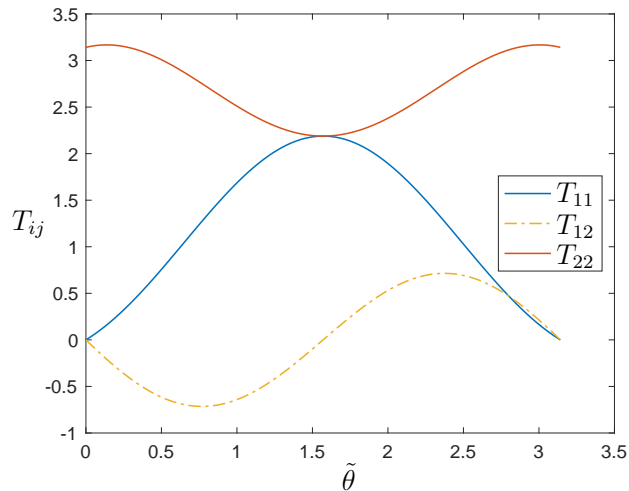}
\caption{Values $T_{11}$, $T_{12}$, and $T_{22}$ in \eqref{matrix_T} as a function of $\tt$. 
}
\label{fig:Tdiag_Txy}
\end{center}
\end{figure}

Finally, combining \eqref{int_theta}, \eqref{Ix} and \eqref{int_J2}, we find that the leading-order contribution to the collision integral is
\begin{align}\label{I_final}
\begin{split}
I &= -\epsilon^ 2 \! \left (D_T  \int_{0}^{\pi}  J \, \d \tt + D_R I_{\tt} \right )= \epsilon^2 \nabla_{\txio} \cdot  \int_0^\pi  \! D \left[\sin \tt{\bf A} + M(\tto,\tt) {\bf B}, \sin \tt  \partial_{\tt}(p p^+) \right] \d \tt.
\end{split}
\end{align}

\subsection{A nonlinear nonlocal diffusion equation} \label{sec:final_eqn}

Inserting the collision integral \eqref{I_final} into \eqref{eq:evolution_eqn_with_bdry_integral}, we find that the integrated Fokker--Planck equation for $N=2$ is
\begin{align}
	\label{eq:fp_final}
	\begin{split}
	\partial_t p =  \nabla_{\xi_1} \cdot \bigg \{ & D \nabla_{\xi_1} p - f(\xi_1) p   + \epsilon ^2 \int_{0}^{\pi} D \left[\sin \theta{\bf A} + M(\theta_1,\theta) {\bf B}, \sin \theta \,  \partial_{\theta} (p p^+) \right] \d \theta \bigg \} .
		\end{split}
\end{align}
The extension from two to $N$ needles is straightforward up to $O(\epsilon^2)$ since only pairwise interactions need to be considered at this order. Noting that the first needle has $(N-1)$ inner regions, one with each of the remaining needles, the marginal density for $N$ needles satisfies
\begin{subequations}
	\label{eq:final_model}
\begin{align}
	\label{eq:fp_final2}
	\partial_t p =  \nabla_{\xi_1} \cdot \left \{ D \nabla_{\xi_1} p - f(\xi_1) p  + \epsilon ^2 (N-1) D \int_0^\pi Q(\theta, p) \d \theta  \right \},
\end{align}
where $D = \diag(D_T, D_T, D_R)$, $f(\xi_1) = (\f_T(\xi_1), f_R(\xi_1)$, and $Q = ({\bf Q}_T,Q_R)$ is given by
\begin{align}
\label{eq:CollisionIntegralContribution}
{\bf Q}_T(\theta, p, p^+)  = \sin \theta{\bf A} + M(\theta_1,\theta) {\bf B}, \quad 
	Q_R(\theta, p, p^+) = \sin \theta p  \partial_{\theta}  p^+.
\end{align}
In \eqref{eq:final_model}, $p = p(\x_1, \theta_1, t)$, $p^+ = p(\x_1, \theta_1 + \theta, t)$, and $M(\theta_1,\theta) =  \Rt T(\theta) \RT$, where $\Rt$ is the rotation matrix by $\theta_1$ (see \eqref{rotation_matrix}) and $T(\theta)$ is the matrix defined in \eqref{matrix_T} (see also Fig.~\ref{fig:Tdiag_Txy}), and 
\begin{align}
\label{AB_final}
\begin{aligned}
	{\bf A} &= \frac{1}{2} \left[  \nabla_{\x_1} (p p^+) + \frac{(p p^+)}{D_T} (\f_T^+ - \f_T)   \right],\\
	{\bf B} &= \frac{1}{2} \left[  p \nabla_{\x_1} p^+ - p^+ \nabla_{\x_1} p + \frac{(p p^+)}{D_T} (\f_T - \f_T^+)   \right].
	\end{aligned}
\end{align}
\end{subequations}
The nonlinearities in \eqref{eq:final_model} encompass the effect that the non-overlap constraint between needles has on the macroscopic dynamics. In particular, we note that the interactions are local in space but nonlocal in angle. 
The integrands ${\bf Q}_T$ and $Q_R$ vanish for $\theta = 0$ (as two parallel needles exclude no volume in phase space), while for $\theta \in (0,\pi)$ they include a series of quadratic terms involving $p$, $p^+$ and their derivatives. The interaction in orientation is of mean-field type (see $I_R$),
where only the ``cross-diffusion'' term $p \partial_{\theta_1} p^+$ appears, whereas in space we obtain full cross-diffusion terms $p \nabla_{\x_1} p^+$ and $p^+ \nabla_{\x_1} p$ as well as a drift-difference term (see ${\bf I}_T$), as in the case of mixtures of hard spheres \cite{Bruna:2012wu}. To give some intuition on their role, consider the kernel ${\bf Q}_T$ for $\theta = \pi/2$ (perpendicular needles). This is the only value for which $T$ is a multiple of the identity (see \eqref{matrix_T}), $T(\pi/2) = \mu I_2$ with $\mu \approx 2.18$. Thus $M(\theta_1, \theta) {\bf B} \equiv  \mu {\bf B}$ and the integrand simplifies to 
\begin{align*}
    {\bf Q}_T(\pi/2, p, p^+) = {\bf A} + \mu {\bf B} = \frac{1}{2}\left[ (\mu + 1) p \nabla_{\x_1} p^+  - (\mu-1) \left(p^+ \nabla_{\x_1} p - \frac{\f_T - \f_T^+}{D_T} \,p p^+ \right) \right].
\end{align*}
In this form, one may readily compare it with the nonlinear terms arising from the interactions between two types of hard-sphere particles of diameter $\epsilon$ (cf. Eq.(22) in \cite{Bruna:2012wu})
\begin{align*}
	{\bf Q}_T(p, p^+) =  \frac{\pi}{2} \left[ 3 p \nabla_{\x_1} p^+  -  p^+ \nabla_{\x_1} p +  \frac{\f_T - \f_T^+}{D_T}  \, p p^+ \right].
\end{align*}
Thus we observe the same structure with an ``effective drift'' $p \nabla_{\x_1} p^+$ due to gradients of the other species, a reduced diffusion $p^+ \nabla_{\x_1} p$ due to concentrations of the other species, and a quadratic drift adjustment with the same relative strength and sign in both needles and hard-spheres cases. The size of the coefficients is larger for hard spheres ($3 \pi/2$ and $\pi/2$) than for needles ($(\mu\pm 1)/2$), as expected given their excluded volume in this specific needles configuration ($\pi$ vs 1).

\begin{remark}[Active Brownian needles] We note that our model \eqref{eq:final_model} may be used to describe a system of $N$ active needles similar to that considered in \cite{MKFL2020} (except that they use the $\theta$-dependent diffusion tensor $\hat D$ in \eqref{orient_tensor}).  
In particular, consider $f_R = 0$ and $f_T (\x, \theta) = v_0 \bf{e}(\theta)$ with $\bf{e}(\theta) = (\cos \theta, \sin \theta)$  in \eqref{sdes}, such that needles drift along their orientation $\theta$ at constant velocity $v_0$. This implies that we must now distinguish between a needle's head and tail as its orientation $\theta$ determines the direction of the drift in position; that is, we must extend the range of $\theta$ to $[0, 2\pi)$. Since the excluded volume between two needles is invariant under switching heads and tails, the terms in \eqref{eq:final_model} that describe the excluded volume, namely $\sin\theta$ and $M(\theta_1, \theta)$ in \eqref{eq:CollisionIntegralContribution} must be extended to $[0, 2\pi)$ as $|\sin \theta|$ and $\tilde M(\theta_1, \theta)$ respectively, where $\tilde M(\theta_1, \theta) = M(\theta_1, \theta)$ for $\theta\in[0, \pi)$ and $\tilde M(\theta_1, \theta) = M(\theta_1, \theta-\pi)$ for $\theta \in[\pi, 2\pi)$. 
Then \eqref{eq:final_model} becomes
\begin{align}
	\label{active_needles}
	\partial_t p =  \nabla_{\x_1} \cdot \left [ D_T \nabla_{\x_1} p - v_0 \e(\theta_1) p  + \phi \int_0^{2\pi} \tilde Q_T(\theta, p) \d \theta  \right ] + D_R \partial_{\theta_1}  \left [ \partial_{\theta_1} p   +  \phi p  \int_0^{2\pi}  \tilde Q_R(\theta, p) \d \theta  \right ] ,
\end{align}
where $\phi = (N-1)\epsilon^2$, with
\begin{align*}
	\tilde Q_T(\theta, p) &= \frac{D_T}{2} \left[ (\tilde M + |\sin \theta|) p \nabla p^+ - (\tilde M - |\sin \theta|) p^+ \nabla p \right]  + \frac{v_0}{2} (\tilde M - |\sin \theta|) p p^+ \widehat \e \\
	&= D_T  \left(\mu^+ p \nabla p^+ - \mu^- p^+ \nabla p \right)   + v_0  \mu^- p p^+ \widehat \e,\\
	\tilde Q_R(\theta,p) &=  |\sin\theta| \partial_\theta p^+
	\end{align*}
where $\widehat \e = \e(\theta_1) - \e(\theta_1 + \theta) $ and $\mu^\pm(\theta_1, \theta) = \frac{1}{2} \left(\tilde M(\theta_1,\theta) \pm |\sin\theta| \right)$.
Rearranging \eqref{active_needles} may be cast in a more familiar form in the active matter community (compare with Eq. (2.29) in \cite{Bruna:2021tb}, corresponding to active Brownian hard disks)
\begin{multline}
	\label{eq_ABN}
\partial_t p + v_0 \nabla \cdot \left[  p (1-\phi  \rho^-) \e(\theta_1) + \phi {\bf m}^- p\right] = D_T \nabla \cdot \left[ (1- \phi \rho^-) \nabla p +  \phi p \nabla \rho^+ \right] \\ + D_R \partial_{\theta_1} \left(\partial_{\theta_1} p + \phi p \bar \rho \right),
\end{multline}
with ``effective'' spatial densities $\rho^\pm$, $\bar \rho$ and magnetisation (also known as polarisation)
\begin{align*}
\rho^\pm = \int_0^{2\pi} \mu^\pm p^+ \d \theta,\qquad \bar \rho = \int_0^{2\pi} \partial_\theta |\sin \theta| p^+ \d \theta, \qquad
{\bf m}^-  = \int_0^{2\pi} \mu^- p^+ {\bf e}^+ \d \theta.
\end{align*}
If the excluded volume between two needles was a constant, then $\rho^\pm \equiv \rho = \int_0^{2\pi} p \d \theta$ (the spatial density), ${\bf m}^- \equiv {\bf m}  = \int_0^{2\pi} p \e \d \theta$ and the nonlinear flux in orientation ($p \bar \rho$) would drop. This is because this term represents changes in orientation brought by the change in excluded volume with relative orientation.
\end{remark}

\section{High rotational diffusion limit} \label{sec:high_rotation}

In the context of colloidal suspensions, the diffusion coefficients corresponding to the rotational and translational motions (parallel or perpendicular to the needle's axis) are not independent. In particular, using Stokes' law, we have that \cite{doi1986theory,Leitmann:2017eb}
\begin{equation} \label{diffusions_coupled}
D_R = 12 D_\perp/\epsilon^2 \qquad D_\parallel = 2 D_\perp,
\end{equation}
where $D_\perp$ and $D_\parallel$ are the translational diffusion coefficients for perpendicular and parallel motion. This means that, instead of the constant diffusion matrix $D = \diag (D_T, D_T, D_R)$ used in our derivation, we would have
\begin{align} \label{orient_tensor}
\begin{aligned}
		\hat D(\theta_1)  = \left(
	\begin{array}{c|c}
    \Rt & \begin{matrix} 0 \\ 0 \end{matrix} \\ \hline
    \begin{matrix} 0 & 0 \end{matrix} & 1
\end{array}\right)	
	  D = \begin{pmatrix}
		\cos \theta_1 D_\parallel & -\sin \theta_1 D_\perp & 0 \\
		\sin \theta_1 D_\parallel & \cos \theta_1 D_\perp& 0 \\
		0 & 0 & D_R
	\end{pmatrix}.
	\end{aligned}
\end{align}
Our derivation can be adapted to allow for a diffusion tensor of this form, resulting in a modified equation for $p$ (in particular, the ${\bf Q}_T$ in \eqref{eq:CollisionIntegralContribution} would change). 
We omit this generalisation here but comment on the asymptotic regime of \eqref{diffusions_coupled}, namely when the rotational diffusion is much larger than the translational diffusion
\begin{align} \label{large_rotation}
	D_R = D_T/\epsilon^2,
\end{align}
and set $D_T \equiv 1$ in this section. 
Inserting \eqref{large_rotation} into \eqref{eq:fp_final2}, we have
\begin{align}
	\label{fp_highDr}
	\begin{split}
		\epsilon ^2 \partial_t p = \ & \epsilon^2  \nabla_{\x_1} \cdot \left[  \nabla_{\x_1} p - \f_T(\xi_1) p \right] + \partial_{\theta_1} \left[ \partial_{\theta_1} p - \epsilon^2 f_R(\xi_1) p \right] \\
		& + \epsilon ^4 (N-1) \nabla_{\x_1} \cdot \int_0^\pi{\bf Q}_T(\theta, p, p^+) \d \theta + \epsilon^2 (N-1) \partial_{\theta_1} \int_0^\pi Q_R(\theta, p, p^+) \d \theta. 
	\end{split}
\end{align}
We look for a solution of \eqref{fp_highDr} of the form $p \sim p_0 + \epsilon^2 p_1 + \cdots$. The leading-order problem gives that $p_0 = p_0 (\x_1,t)$, that is, the leading-order problem is independent of angle. Collecting the $O(\epsilon^2)$-terms in \eqref{fp_highDr} yields
\begin{align}
	\label{fp_highDr2}
\partial_t p_0 =   \nabla_{\x_1} \cdot \left[  \nabla_{\x_1} p_0 - \f_T(\xi_1) p_0 \right] + \partial_{\theta_1} \left[ \partial_{\theta_1} p_1 - f_R(\xi_1) p_0 \right], 
\end{align}
where we have used that $Q_R(\theta, p_0, p_0^+) \equiv 0$. The $O(\epsilon^4)$ of \eqref{fp_highDr} is 
\begin{align}
	\label{fp_highDr4}
	\begin{split}
		\partial_t p_1 = \ &  \nabla_{\x_1} \cdot \left[  \nabla_{\x_1} p_1 - \f_T(\xi_1) p_1 \right] + \partial_{\theta_1} \left[ \partial_{\theta_1} p_2 -  f_R(\xi_1) p_1 \right] \\
		& +  (N-1) \nabla_{\x_1} \cdot \int_0^\pi {\bf Q}_T(\theta, p_0, p_0^+)\d \theta +  (N-1) \partial_{\theta_1}\int_0^\pi  Q_R(\theta, p_0, p_1^+) \d \theta,
	\end{split}
\end{align}
noting that $Q_R(\theta, p_1, p_0^+) \equiv 0$. We now write an equation for the spatial density
\begin{equation*}
	\rho (\x_1, t) := \int_0 ^\pi ( p_0 + \epsilon^2 p_1) \d \theta_1.
\end{equation*}
Combining \eqref{fp_highDr2} and \eqref{fp_highDr4}, and using periodicity in $\theta_1$, we find
\begin{align}
	\label{fp_highDr_final}
		\partial_t \rho = \nabla_{\x_1} \cdot  \left[  \nabla_{\x_1} \rho -  \! \int_0^\pi \f_T(\xi_1) p(\xi_1,t) \d \theta_1 + \epsilon^2 (N-1) \!\int_0^\pi \int_0^\pi {\bf Q}_T(\theta; p_0, p_0^+) \d \theta \d \theta_1  \right].
\end{align}
In particular, if we assume that $\f_T$ is independent of angle, then 
\begin{align}
	\label{fp_highDr_final2}
		\partial_t \rho = \nabla_{\x_1} \cdot  \left[  \nabla_{\x_1} \rho -  \f_T(\x_1) \rho + \epsilon^2 (N-1) \!\int_0^\pi \int_0^\pi  {\bf Q}_T(\theta; p_0, p_0^+) \d \theta \d \theta_1  \right].
\end{align}
Using that $p_0 = p_0^+$ and $\f_T = \f_T^+$, from \eqref{AB_final} we have that ${\bf A} = p \nabla_{\x_1} p$ and ${\bf B}=0$ and hence ${\bf Q}_T = \sin \theta \nabla_{\x_1} (p_0^2)/2$. The double integral on ${\bf Q}_T$ is then $\pi \nabla_{\x_1} (p_0^2)\sim \frac{1}{\pi} \nabla_{\x_1} (\rho^2)$ using that $\rho = \pi p_0 + O(\epsilon^2)$. We find that  \eqref{fp_highDr_final2} reduces to
\begin{align}
	\label{fp_highDr_final3}
		\partial_t \rho = \nabla_{\x_1} \cdot  \left\{ \left[1  + \frac{2}{\pi} (N-1) \epsilon^2 \rho \right] \nabla_{\x_1}\rho -  \f_T(\x_1) \rho  \right \}.
\end{align}
Therefore we find that the equation satisfied by $N$ needles of length $\epsilon$ in the limit of large rotational diffusion is a nonlinear diffusion equation of the same form as the equation \eqref{al:brunachapman_equation} satisfied by $N$ disks of diameter $\epsilon$. 
 Comparing the two equations, we have that the effective diameter of a needle with very fast rotational diffusion is $\sqrt{2}/\pi$ times its length $\epsilon$. That is, the needle excludes roughly 45\% less volume than a disk of diameter $\epsilon$.

\section{Space homogeneous solutions} \label{sec:space_hom}

In this section, we consider spatially homogeneous solutions to \eqref{eq:final_model}, that is,  solutions of the form $p(\xi_1,t) =  p(\theta_1,t)$ satisfying
\begin{align} \label{spatially_hom}
	D_R^{-1}\partial_t p = \partial_{\theta_1}^2 p + \epsilon ^2 (N-1) \partial_{\theta_1} \left(p \int_{0}^{\pi} \sin \theta \partial_{\theta} p^+ \d \theta \right).
\end{align}
The integral is
\begin{align*}
	\int_{0}^{\pi} \sin \theta \partial_{\theta} p^+ \d \theta = - \int_{0}^{\pi} \cos \theta p(\theta_1 + \theta) \d \theta = \int_{0}^{\pi} \cos \theta p(\theta_1 - \theta) \d \theta =  W' \ast p,  
\end{align*}
where $W(\theta) = \sin(\theta)$. Therefore,  the space-homogeneous system of interacting needles of length $\epsilon$ is described by a periodic  McKean--Vlasov equation with an attractive potential $W$ (see, e.g., \cite{carrillo2020long,kuramoto1981rhythms})
\begin{align} \label{MKV}
	D_R^{-1}\partial_t p = \partial_{\theta_1}^2 p + \epsilon ^2 (N-1) \partial_{\theta_1} (p W' \ast p).
\end{align}
We study the linear stability of the homogeneous solution $p_* = 1/\pi$ of \eqref{MKV} by considering a perturbation of the form
$$
p = p_* + \delta e^{\lambda t} \sum_{n \ge 0} a_n \cos(2n\theta_1) + b_n \sin(2n \theta_1),
$$ 
with $\delta \ll 1$. Inserting this into \eqref{MKV}, linearising and keeping terms of $O(\delta)$, we arrive at
\begin{align*}
	\lambda  = -4n^2 D_R \left(1 - \frac{2 \phi n}{(4n^2-1) \pi} \right),
\end{align*}
where $\phi = \epsilon ^2 (N-1)$. We look for growing modes by imposing  $\lambda>0$, leading to $2\phi n > (4 n^2-1)\pi$. The most unstable mode ($n=1$) leads to linear instability if 
\begin{equation}
	\label{critical_value}
\phi > \phi_c = \frac{3 \pi}{2}.
\end{equation}
Note that, while $\phi$ represents an effective volume fraction (which would be bounded for isotropic bodies by their close packing densities, e.g., $\phi <0.74$ for closely packed hard disks in two dimensions), the hard-core needle system admits any $\phi\in[0, \infty)$, with $\infty$ corresponding to a system of perfectly aligned needles.

It is also worth pointing out that, while our derivation relied on a diluteness assumption $\phi\ll 1$, the critical volume fraction is $\phi_c = O(1)$. Therefore, the aggregation behaviour occurs outside the region of validity of our PDE model \eqref{eq:final_model} and, as a by-product, of the space-homogeneous model \eqref{MKV}. In fact, the value $\phi_c$ agrees with the bifurcation point of isotropic-nematic transition obtained in \cite{PhysRevA.17.2067} using Onsager's theory of orientational order \cite{Onsager:1949jk}. In particular, Onsager considers the virial expansion of the orientational probability density up to the second virial coefficient (which depends on two-particle interactions, and Onsager obtains for a variety of hard anisotropic particles evaluating the excluded volume for a pair of such particles). While the third- and higher-order virial coefficients are negligible for hard needles in $\mathbb R^3$, it is not the case in the present case of two dimensions \cite{PhysRevA.17.2067}. Therefore, the value we obtain for $\phi_c$ should be taken with caution, and indeed Monte Carlo simulations have found the critical density at the transition to be $\phi_c\approx 7$ \cite{frenkel1985evidence}.

The stationary solutions of \eqref{MKV} satisfy
$$
\partial_{\theta_1} p_s + \phi p_s W' \ast p_s = -J,
$$
where $J$ is a constant corresponding to the flux of the stationary solution. Without any external forcing, we expect solutions with $J=0$. Imposing $J=0$ and integrating, we arrive at
\begin{equation*}
p_s(\theta_1) = C \exp \left( \phi \int_0^{\theta_1} (W' \ast p) (\theta) \d \theta   \right),
\end{equation*}
where $C$ is a normalisation constant such that $\int_0^\pi p_s \d \theta_1 = 1$. We consider a fixed-point iteration method to compute $p_s(\theta_1)$ above for various values of $\phi$. Specifically, given an initialisation $p_0$ (normalised to one), we compute 
\begin{equation} \label{fpm}
	p_{k+1} = C \exp \left( \phi \int_0^{\theta_1} (W' \ast p_k) (\theta)\, \d \theta   \right), \qquad \text{for } k = 1, 2, \dots.
\end{equation}
We initialise the scheme with the most unstable mode from the linear stability analysis ($p_0(\theta_1) = 1/\pi + \delta \cos(2 \theta_1)$) and solve \eqref{fpm} with Chebfun \cite{Driscoll2014} until it reaches a stationary profile. We consider several values of $\phi\ge 3\pi/2$ so that we expect nontrivial stationary states.  Figure \ref{fig:stationary_sols} shows the results for ten values of $\phi$. We observe that the stationary solution becomes more concentrated as $\phi$ increases. This means that needles are forced to align more to avoid overlapping as their number increases.
\def \scl {1.0}
\begin{figure}[htb]
    \begin{center}
    \psfrag{p}[][t][\scl][-90]{$p_s$}
    \psfrag{pt}[][t][\scl][-90]{$p$}
    \psfrag{t}[][b][\scl]{$\theta_1$}
    \psfrag{f}[][b][\scl]{$\phi$}
    \psfrag{a}[][][\scl]{\ (a)}
    \psfrag{b}[][][\scl]{\ (b)}
    \psfrag{d}[][b][\scl]{$t$}
    \includegraphics[width = 0.475\textwidth]{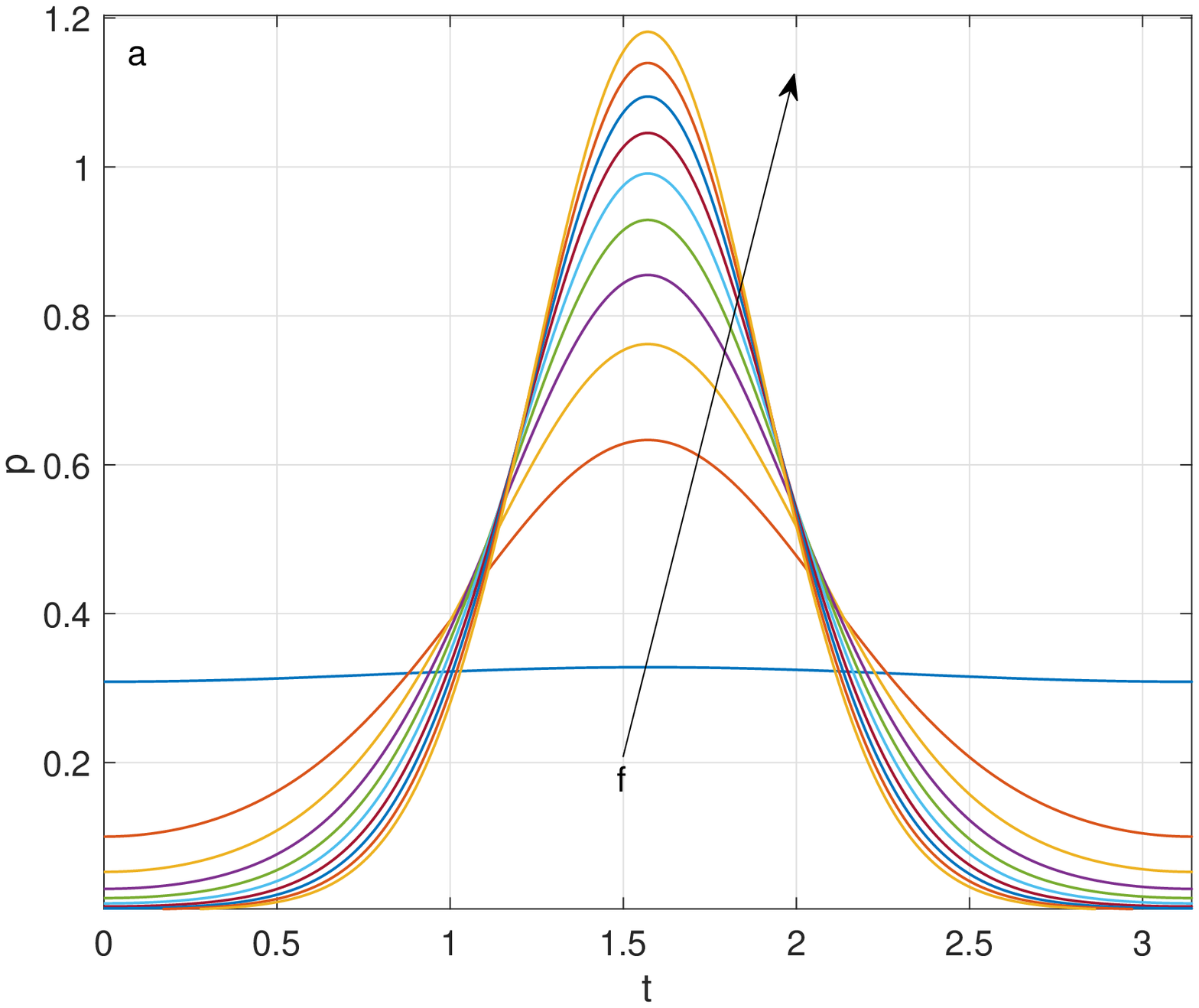}
    \quad
    \includegraphics[width = 0.475\textwidth]{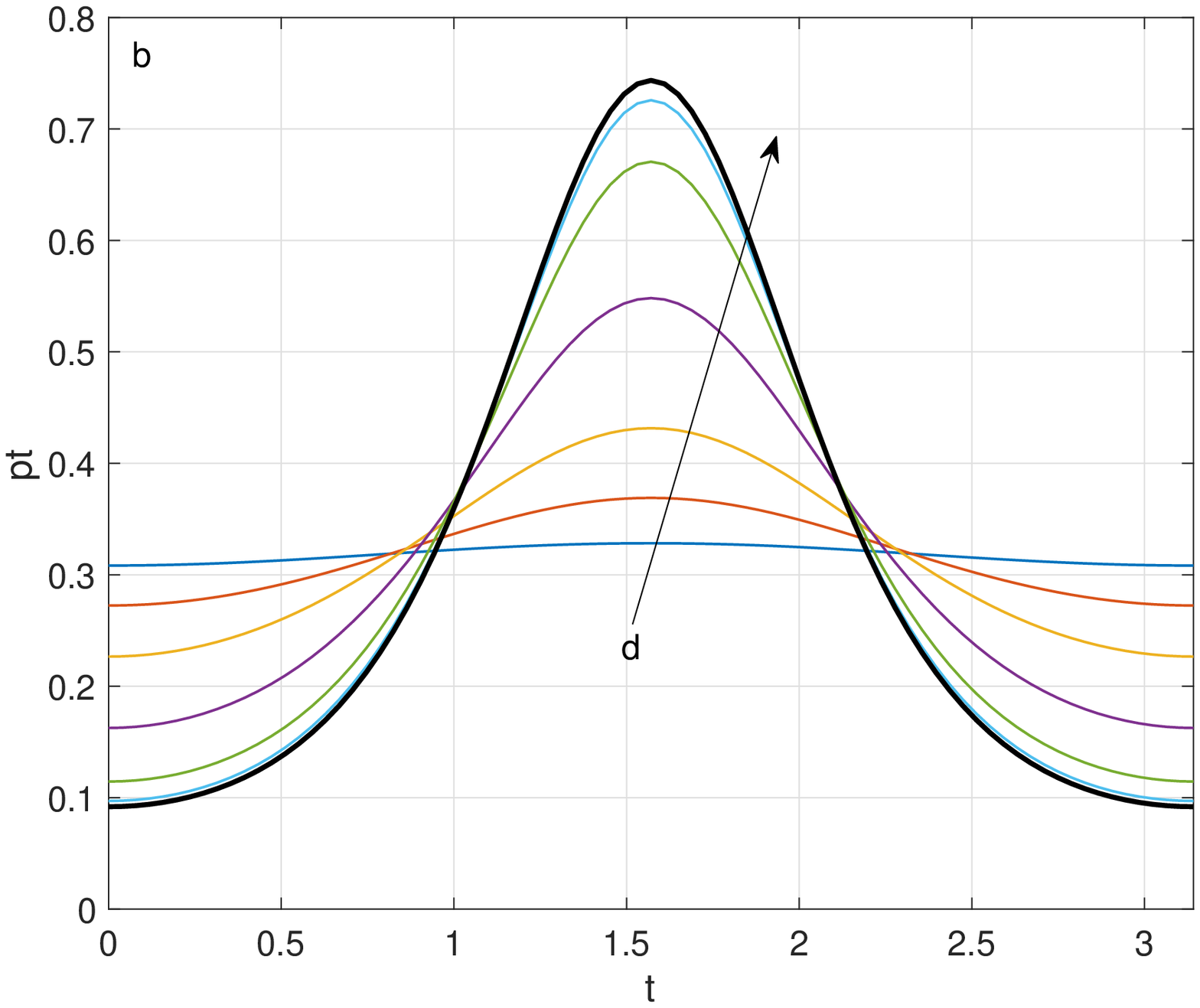}
    \end{center}
  \caption{(a) Stationary solutions $p_s(\theta_1)$ of the space-homogeneous problem \eqref{MKV} for different values of $\phi = (3\pi/2) + k/2$ for $k = 0, \dots, 10$. Solutions are obtained via a fixed-point iterative scheme using Chebfun \cite{Driscoll2014}. (b) Time evolution $p(\theta_1,t)$ for $\phi = 1.1 \times 3\pi/2$ for a small initial perturbation $p_0 = \pi^{-1} - 0.01\cos(2\theta_1)$. Times shown are $t= 0, 4, 6, 8, 10, 12, 20$. At $t=20$, the solution has already reached the stable equilibrium.}
 \label{fig:stationary_sols}
\end{figure}

\section{Discussion}
We have systematically derived an effective PDE model for a system of non-overlapping Brownian needles in two dimensions \eqref{eq:final_model}. The nonlinearities of the PDE describe the effect of pairwise interactions at the macroscopic level: interactions are nonlocal in angle (the nonlinearity is of mean-field type, only $p \partial_{\theta_1} p^+$ term) and local in position (full cross-diffusion terms $p \nabla_{\x_1} p^+$ and $p^+ \nabla_{\x_1} p$ as well as a drift-difference term appear, consistent with other local-in-space models \cite{Bruna:2012wu,Mason.2022}). To gain insight into the behaviour of the PDE model, we consider two simplifications. First, we obtain a reduced PDE for the spatial density in the high-rotational diffusion limit. By comparing the resulting PDE with the effective PDE for hard-core disks in two dimensions, we find that the needles' effective diameter is about 45 per cent of their length. Second, we consider space-homogeneous solutions of the nonlocal PDE and show they satisfy a well-known McKean-Vlasov equation with an attractive potential in orientation. Notably, we identify an instability of the uniform distribution (in angle) for effective packing densities above a critical threshold, see \eqref{critical_value}. Intuitively, we expect this phase transition to occur and arise from the finite-size interactions between needles. Indeed, the instability corresponds to the emergence of a preferred direction of needles to exclude less volume in configuration space in crowded settings. 
Let us point out that the nonlocal interaction term in Eq. \eqref{MKV} includes the size of the excluded volume. 

In this work, we find that the strength of the nonlinearity in the macroscopic PDE is proportional to the total excluded region volume $(N-1)\epsilon^2 \sin\theta$. The form of such nonlinearity is nontrivial in the full PDE \eqref{eq:final_model} (due to the spatial interactions). Still, it may have been inferred in the space-homogeneous case \eqref{spatially_hom} (in fact, this was the approach taken in \cite{PhysRevA.17.2067} using Onsager's free energy functional based on the geometry of the excluded region). 
A natural question is whether this can be generalised to similar systems. A particularly interesting case is that of Brownian needles in three dimensions, which have zero excluded volume in configuration space (for fixed relative angles, the excluded region is a two-dimensional surface in $\R^3$). If the result from two dimensions were to extend to three dimensions, it would imply that the effective PDE for needles in three dimensions would not ``see'' the non-overlapping constraints, at least not to $O(N\epsilon^3)$.







\vspace{1cm}

\textbf{Acknowledgements} M. Bruna was supported by a Royal Society University Research 
Fellowship (grant no. URF/R1/180040). The authors would like to thank 
Martin Burger for the helpful discussions.



\appendix
\section{Solution of the first-order inner problem via conformal mapping} \label{sec:conformal}

We solve problems \eqref{prob_u1} and \eqref{prob_u2} by mapping them to problems in the interior of a circle. 
We consider the problem for $u_1$ \eqref{prob_u1}; the problem for $u_2$ follows similarly.

Let $D$ denote the exterior of the rhombus in  the $z$-plane, $D = \mathbb C \setminus \rhombus$, where $z = \tilde x + \im \tilde y$, and let $\Gamma$ be its boundary, $\Gamma = \partial \rhombus$. Let $\Delta_z$ denote the Laplacian operator $\partial^2/\partial \tilde x^2 + \partial^2/\partial \tilde y^2$.  
We look for a complex function $w_1: D \to D$ such that the solution we need is given as $u_1 = \text{Re}(w_1)$. 
By the Cauchy--Riemann relations it follows that the boundary condition $\nabla_{\tx} u_1 \cdot \tilde \n =  0$ on $\Gamma$ is equivalent to imposing that the conjugate harmonic function $\text{Im}(w_1)$ is constant on $\Gamma$, for example, equal to zero. Then $w_1$ must satisfy
\begin{equation} \label{prob_w1}
	\begin{alignedat}{2}
\Delta_z w_1	 & = 0  & & \text{in } D, \\
\text{Im}(w_1) &=  0   &\qquad & \text{on } \Gamma,\\
w_1 & \sim z & &  \text{at } \infty.
\end{alignedat}
\end{equation}

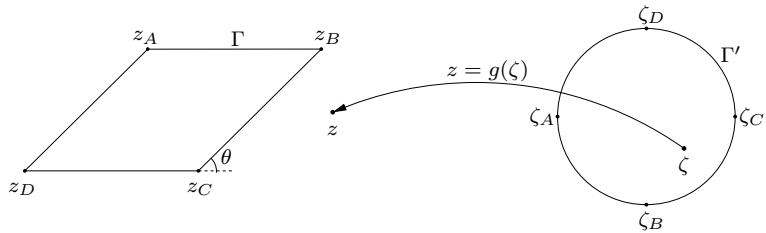
\begin{figure}
\begin{center}
\input{conformal.pstex_t} 
\caption{Mapping of the exterior of the unit circle with boundary $\Gamma'$ into the D of the rhombus with boundary $\Gamma$. }
\label{fig:conformal}
\end{center}
\end{figure}

To proceed with the solution of \eqref{prob_w1}, we seek a transformation that simplifies the definition domain. In particular, we look for an analytic function $z = g(\zeta)$ that maps a domain $D'$ of the $\zeta$ plane, namely the interior of the unit disk, to $D$ in the $z$ plane (see Fig.~\ref{fig:conformal}). Then the unit circle, denoted by $\Gamma'$, is mapped into the boundary of the rhombus $\Gamma$. 
 This is a Schwarz--Christoffel transformation, given by \cite[eq. (4.6)]{driscoll2002schwarz}
\begin{equation}
	\label{g_transform}
	z = g(\zeta) = a_0 + a(\tt) \int^\zeta (1- t^2)^{\tt/\pi} (1 + t^2)^{1-\tt/\pi} t^{-2} \, \d t,
\end{equation} 
where $a_0$ and $a(\tt)$ are chosen so that $g(\zeta_k) = z_k$, for $k = A, B, C, D$, where $\zeta_k  = \pm 1, \pm \im$ (see \eqref{fig:conformal}). Note that, as we move through the points $A\to B\to C\to D\to A$, we travel the circle counterclockwise but the rhombus clockwise (so that both curves are positively oriented, i.e., we have the domain to our left as we travel on its boundary). 
We note that $g(\zeta)$ goes to infinity like $-a(\tt)/\zeta$ as $\zeta \to 0$.  The constant $a(\tt)$ is given exactly as
\begin{equation} \label{a_eq}
	a(\tt)=  \frac{\alpha}{\beta - \im \gamma},
\end{equation}
where $\alpha$, $\beta$ and $\gamma$ are the following real functions of $\tt$:
\begin{align*}
	\alpha(\tt) &=  2^{1+ 2\tt/\pi} \sec{\tt},\\
	\beta(\tt) &= \Gamma\left( \tfrac{1}{2} - \tfrac{\tt}{\pi} \right)\Gamma\left( 1+\tfrac{2\tt}{\pi} \right)\\
	& \quad \times\left[ {}_2 \bar F_1\left( \frac{1}{2}, \frac{\tt}{\pi};\frac{3}{2} + \frac{\tt}{\pi};-1 \right)- 2 {}_2 \bar F_1\left( -\frac{1}{2}, \frac{\tt}{\pi};\frac{1}{2}+ \frac{\tt}{\pi};-1 \right)  \right],\\
	\gamma(\tt) &= 16^{\tt/\pi}\Gamma\left( \tfrac{1}{2} + \tfrac{\tt}{\pi} \right)\Gamma\left(1- \tfrac{2\tt}{\pi} \right)\\
	& \quad \times\left[ {}_2 \bar F_1\left( \frac{1}{2}, -\frac{\tt}{\pi}; \frac{3}{2}-\frac{\tt}{\pi};-1 \right)+ 2 {}_2 \bar F_1\left( -\frac{1}{2}, -\frac{\tt}{\pi}; \frac{1}{2}-\frac{\tt}{\pi};-1 \right)  \right],
\end{align*}
where ${}_2 \bar F_1 (a, b;c;z)= {}_2 F_1 (a, b;c;z) / \Gamma(c)$ is the regularised hypergeometric function.

The map $g$ corresponding to $\tt = \pi/4$ is illustrated in \eqref{fig:SC_map}(a), and the complex constant $a(\tt) = a_1(\tt) + \im a_2(\tt)$, where $a_1 = \alpha \beta/(\beta^2 + \gamma^2)$ and $a_2 = \alpha \gamma/(\beta^2 + \gamma^2)$, is shown in  \eqref{fig:SC_map}(b). Note that although $\alpha, \beta, \gamma$ are singular at $\tt = \pi/2$, $a_1$ and $a_2$ are not.
\def \scl {.9}
\begin{figure}
\unitlength=1cm
\begin{center}
\psfrag{a}[][][\scl]{(a)} \psfrag{c}[][][\scl]{(b)} \psfrag{ai}[][][\scl][-90]{$a_i$} \psfrag{a1}[][][\scl]{$a_1$} \psfrag{a2}[][][\scl]{$a_2$}  
\psfrag{th}[][][\scl]{$\tt$} 
\includegraphics[height=.4\textwidth]{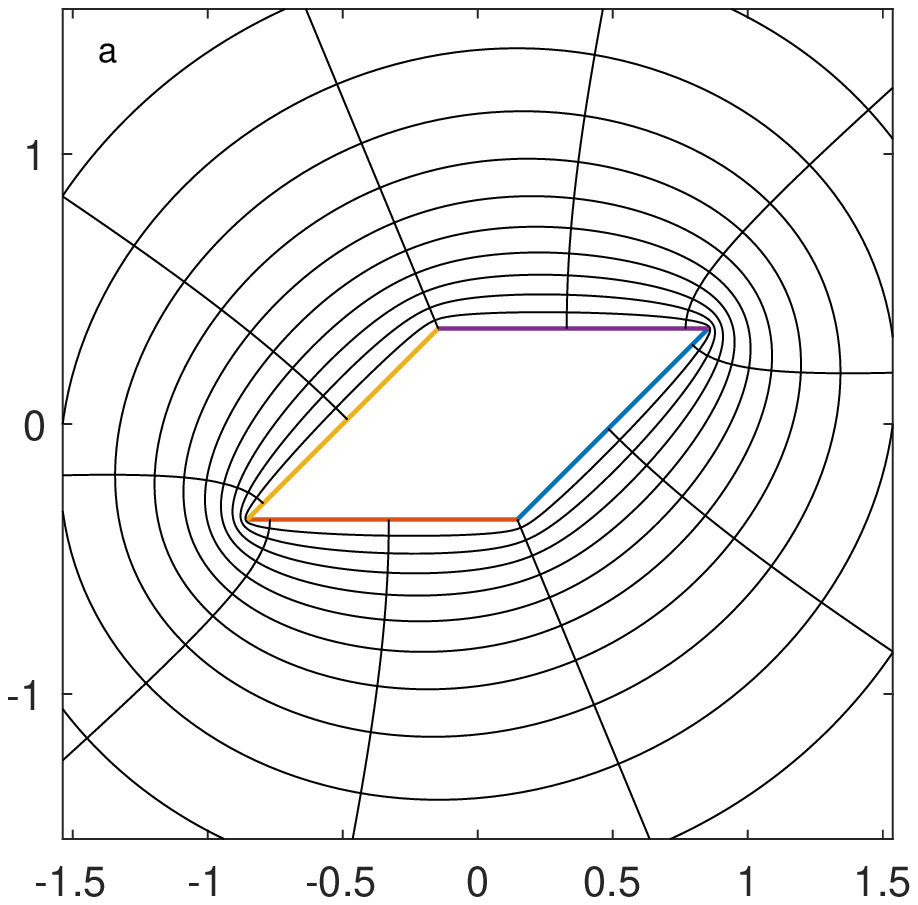} \quad
\includegraphics[height=.4\textwidth]{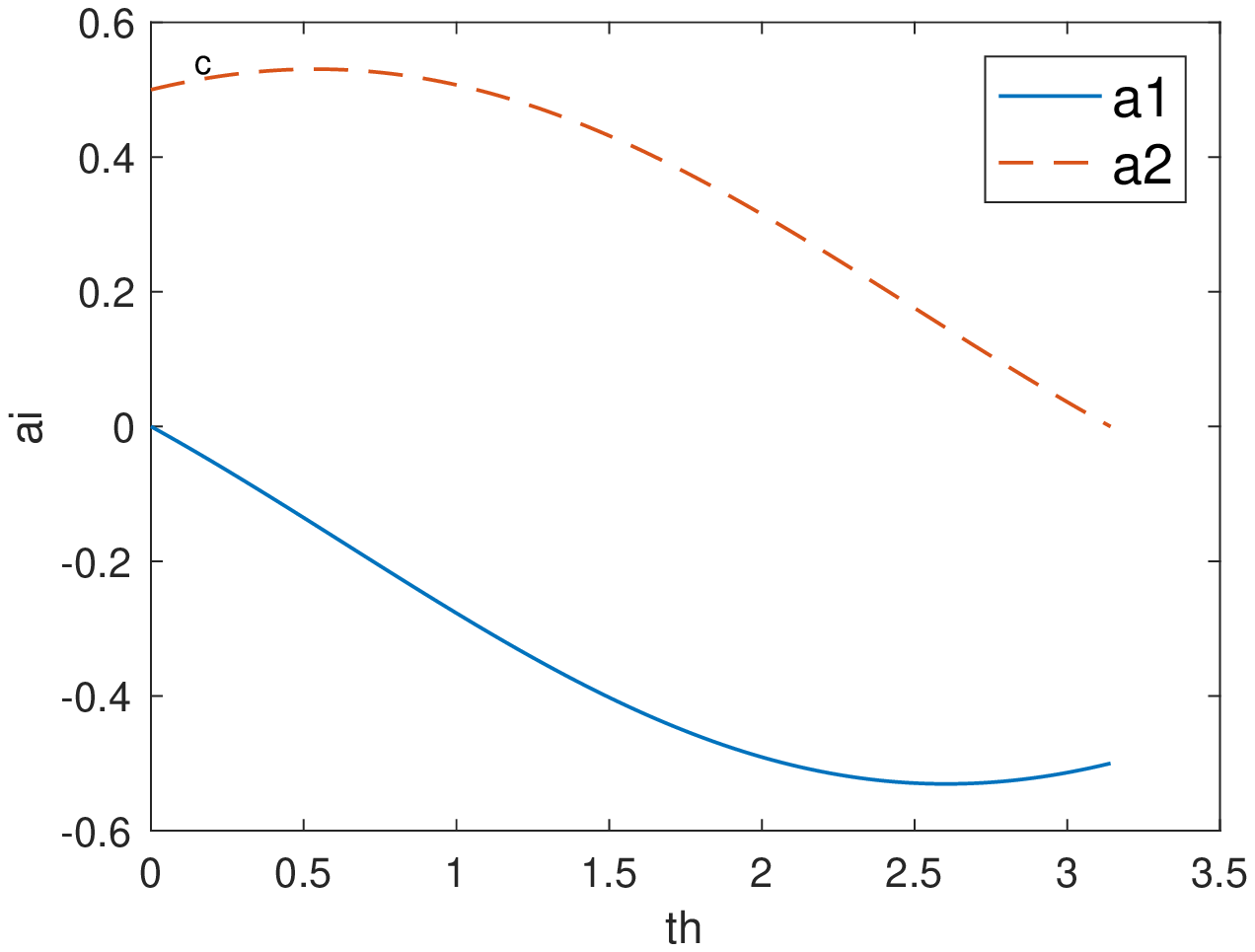} 
\caption{(a) Schwarz--Christoffel map $g$ in \eqref{g_transform} from the interior of unit circle to the exterior of the rhombus, for $\tt = \pi/4$. The black curves are the images of ten evenly spaced circles centred at the origin and ten evenly spaced radii in the unit disk. Plot generated using the Schwarz--Christoffel MATLAB Toolbox \cite{driscoll2005schwarz}. (b) Real and imaginary parts of the multiplicative constant $a(\tt)$ \eqref{a_eq}.}
\label{fig:SC_map}
\end{center}
\end{figure}

We now write the problem in the $\zeta$ plane. If $w_1$ satisfies \eqref{prob_w1} in $D$, $W_1(\zeta):= w_1(g(\zeta))$ is satisfies the following problem in $D'$:
\begin{equation} \label{prob_W1}
	\begin{alignedat}{2}
	\Delta_\zeta W_1 &=0 && |\zeta|<1,\\
\text{Im}(W_1) &=  0   &\qquad & |\zeta|=1,\\
W_1  &\sim -a(\tt) \zeta^{-1} & &  \text{at } 0,
\end{alignedat}
\end{equation}
where $\Delta_\zeta$ denotes the Laplacian operator in the $\zeta$-plane. 
The solution to \eqref{prob_W1} is
\begin{equation}
	\label{W1_sol}
	W_1(\zeta) = -\left(\overline a(\tt)\zeta + \frac{a(\tt)}{\zeta}\right). 
\end{equation}

Repeating the same procedure to solve for \eqref{prob_u2}, we find that $u_2 = \text{Re}(w_2)$ where $w_2$ satisfies \eqref{prob_w1} but replacing the condition at infinity by $w_2 \sim - \im z$. Then the solution in the $\zeta$ plane $W_2(\zeta):= w_2(g(\zeta))$ needs to go like $i a(\tt)/\zeta$ at the origin and is therefore is given by 
\begin{equation}
	\label{W2_sol}
	W_2(\zeta) =  - \im \left( \overline a(\tt) \zeta - \frac{a(\tt)}{\zeta} \right). 
\end{equation}

\section{Collision integral} \label{sec:integral_all}

In this appendix, we evaluate the integral $J$ in \eqref{Ix}, 
\begin{equation}
	J =  \int_{\partial \rhombus}  \RT \nabla_{\txo} \tPo{1} \cdot \tilde \n \, \d S_{\tx},
\end{equation}
where $\rhombus$ is the excluded rhombus in the inner region with $|\rhombus| = \sin \tt$ (see Remark \eqref{rem:excluded_region}). 
Using the first-order inner solution $\tPo{1}$ \eqref{p1_decomposed}, we write 
 $J = J_A + J_B + J_C$ with 
\begin{align} \label{int_J}
\begin{aligned}
J_A &= \int_{\partial \rhombus}  \RT \nabla_{\txo} (\RT {\bf A}\cdot \tx) \cdot \tilde \n \, \d S_{\tx} = \nabla_{\txo} \cdot \left[ \Rt  \left( \int_{\partial \rhombus} \tx \otimes \tilde \n \, \d S_{\tx} \right) \RT  {\bf A} \right] ,\\
J_B &= \int_{\partial \rhombus}  \RT \nabla_{\txo} (\RT {\bf B} \cdot {\bf u}) \cdot \tilde \n \, \d S_{\tx} = \nabla_{\txo} \cdot \left[ \Rt  \left( \int_{\partial \rhombus} {\bf u} \otimes \tilde \n \, \d S_{\tx} \right) \RT  {\bf B} \right],\\
J_C &= \RT \nabla_{\txo} C_\infty  \cdot \int_{\partial \rhombus} \tilde \n \, \d S_{\tx}.\\
	\end{aligned}
\end{align}
We have $J_C = 0$ since we integrate the normal along the closed curve $\partial \rhombus$. To evaluate $J_A$ and $J_B$, we are left to compute the matrices inside the round brackets, which we denote by $-Q$ and $-T$, respectively,
$$
Q = -\int_{\partial \rhombus} \tx \otimes \tilde \n \, \d S_{\tx}, \qquad T = - \int_{\partial \rhombus} {\bf u} \otimes \tilde \n \, \d S_{\tx}.
$$
The first row of $Q$ is
\begin{equation} \label{int_J1}
Q_{1\cdot} = -
\int_{\partial \rhombus} \tilde x \tilde \n \, \d S_{\tx} \sim  \int_{\rhombus} \nabla_\tx \tilde x  \, \d {\tx} =  (1,0) \int_{\rhombus} \d {\tx} =  (1, 0) \sin \tt,
\end{equation}
applying the divergence theorem (on $\tilde x {\bf c}$ with $\bf c$ constant). 
The $\sim$ equivalence is due to the fact that $\tilde \n$ is the projection of the unit normal $\tilde n$ onto the $\tx$ plane, and so it is not normalised (see \eqref{rem:normal_transform} and discussion thereafter). However, since the component of $\tilde n$ in the $\tt$ direction is $O(\epsilon)$, and we only require the leading order of $J$, we can treat $\tilde \n$ as if it were the unit normal on $\rhombus$. For example, on the top edge of the rhombus, we have $\tilde \n \sim (0, -1)$ (Fig.~\ref{fig:excluded_volume}). 
Note also the change in sign in the first equivalence since $\tilde \n$ is the inward unit normal to $\rhombus$. Similarly,  we find that the second row of $Q$ is 
\begin{equation}
	\label{int_J2_app}
	Q_{2\cdot} = -\int_{\partial \rhombus} \tilde y \tilde \n \, \d S_{\tx} \sim (0, 1) \sin \tt.
\end{equation}
Therefore, we find that $Q(\tt) = \sin\tt I_2$.
Matrix $T$ has rows 
\begin{equation}\label{ui_integral}
	T_{i\cdot} = - \int_{\partial \rhombus} u_i (\tx) \tilde \n \, \d S_{\tx},
\end{equation}
where $u_i$ for $i = 1, 2$ solve \eqref{prob_u1} and \eqref{prob_u2}, respectively. 
Rather than transforming the solutions $W_1$ and $W_2$ obtained in Appendix \ref{sec:conformal} back to the $\tx$ plane, we express the integrals as complex integrals in the $\zeta$ plane (see Figure \ref{fig:conformal}). 
To transform \eqref{ui_integral} into a complex integral, first recall that $z = \tilde x + \im  \tilde y$. Given a parameterisation $(\tilde x(s), \tilde y(s))$ of $\partial \rhombus \equiv \Gamma$, the integral along the arclength is 
$$
\d S_\tx \equiv \d {\bf s} = (\tilde x'(s), \tilde y'(s)) \d s = (\tilde x'(s) + \im \tilde y'(s)) \d s = z'(s) \d s = \d z.  
$$
Since the curves $\Gamma$ and $\Gamma'$ are positively oriented (see Figure \ref{fig:conformal}), the corresponding outward normals to the interior of the rhombus or the exterior of the circle, respectively, are given by a $-\pi/2$ rotation, or $- \im$, of the tangent vector. 
Therefore $T_{i}$  as a complex integral is
\begin{equation} \label{integral_J3}
	T_{i \cdot} =  \im \int_{\Gamma} u_i (z) \, \d z =  \im \int_\Gamma w_i (z) \, \d z =  \im \int_{\Gamma'} W_i(\zeta) g'(\zeta) \d \zeta.
\end{equation}
In the second equality we have used that $\text{Im}(w_i) = 0$ on $\Gamma$ (see \eqref{prob_w1}) and in the third that $W_i(\zeta) = w_i(g(\zeta))$.

\begin{figure}[tb]
\begin{center}
\includegraphics[width=.3\textwidth]{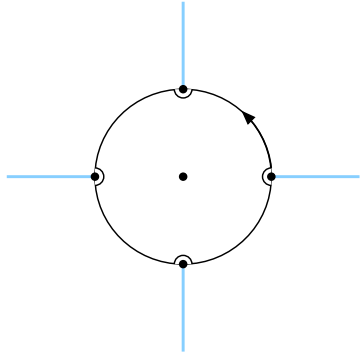} 
\caption{Contour to compute the integral \eqref{integral_J3}.}
\label{fig:contour}
\end{center}
\end{figure}
The integrand in \eqref{integral_J3} has a singularity at the origin and branch points at $\pm 1$ and $\pm \im$. We choose branch cuts going to infinity so that the contour of integration follows $\Gamma'$ with four small semicircular indentations at the branch points as shown in \eqref{fig:contour}. This way, $T_{i\cdot}$ can be computed using Cauchy's Residue Theorem, with $2 \pi \im $ times the residue at the origin and $-\pi \im $ times the residues at the four branch points.\footnote{Note that the contribution of the four points on the unit circle is $-\pi \im $ since it is only half a circle and we are taking the small semicircles clockwise.} In fact, the four branch points do not contribute to the integral for $\tt\in (0,\pi)$ as their residues are zero (no singularities). Because of the form of $W_1(\zeta)$ and $W_2(\zeta)$ (see \eqref{W1_sol} and \eqref{W2_sol}), it is sufficient to compute the following residues
\begin{equation}
	\label{residues}
	\Res_{\zeta = 0} \left[ \zeta g'(\zeta) \right] = a(\tt), \qquad
	\Res_{\zeta = 0} \left[ \zeta^{-1} g'(\zeta) \right] = \left(1 - \frac{2\tt}{\pi}  \right) a(\tt).
\end{equation}
Substituting in the expressions for $W_i$ \eqref{W1_sol} and \eqref{W2_sol} in \eqref{integral_J3} and using  \eqref{residues}, we find
\begin{align} \label{J3_result}
\begin{aligned}
	T_{1\cdot} &=  -\im  \int_{\Gamma'} \left(\overline a(\tt)\zeta + \frac{a(\tt)}{\zeta}\right) g'(\zeta) \d \zeta = 2 \pi \left \{ \overline a \Res_{\zeta = 0} \left[ \zeta g'(\zeta) \right] + a \Res_{\zeta = 0} \left[ \zeta^{-1} g'(\zeta) \right] \right \}\\
	& = 2\pi a \overline a - (4\tt -2 \pi) a^2,
	\end{aligned}
\end{align}
and
\begin{align} \label{J4_result}
\begin{aligned}
	T_{2\cdot} &=  \!\int_{\Gamma'} \left(\overline a(\tt)\zeta - \frac{a(\tt)}{\zeta}\right) g'(\zeta) \d \zeta = 2 \pi \im   \left \{ \overline a \Res_{\zeta = 0} \left[ \zeta g'(\zeta) \right] -  a \Res_{\zeta = 0} \left[ \zeta^{-1} g'(\zeta) \right] \right \}\\
	& = 2\pi \im a \overline a +\im (4\tt -2 \pi) a^2. 
	\end{aligned}
\end{align}
Writing \eqref{J3_result} and \eqref{J4_result} as two-dimensional vectors, we obtain the symmetric matrix
\begin{equation} \label{matrix_T}
	T(\tt) := - \int_{\partial \rhombus} {\bf u} \otimes \tilde \n \, \d S_{\tx} = 4 \begin{bmatrix}
 	a_1^2(\pi-\tt)  + a_2^2 \tt  & a_1 a_2 (\pi-2  \tt) \\
 	a_1 a_2 (\pi-2  \tt) & \ a_2^2 (\pi-\tt )  + a_1^2 \tt
 \end{bmatrix},
\end{equation}
where $a_1$ and $a_2$ are shown in Figure \ref{fig:SC_map}(b).


\def\cprime{$'$}

\end{document}

%% file: excluded_region.pstex_t
\begin{picture}(0,0)%
\includegraphics{excluded_region.pstex}%
\end{picture}%
\setlength{\unitlength}{2368sp}%
\begingroup\makeatletter\ifx\SetFigFont\undefined%
\gdef\SetFigFont#1#2#3#4#5{%
  \reset@font\fontsize{#1}{#2pt}%
  \fontfamily{#3}\fontseries{#4}\fontshape{#5}%
  \selectfont}%
\fi\endgroup%
\begin{picture}(4752,4104)(-2339,-853)
\put(1801,1814){\makebox(0,0)[lb]{\smash{{\SetFigFont{7}{8.4}{\rmdefault}{\mddefault}{\updefault}{\color[rgb]{0,0,0}$\theta$}%
}}}}
\put(151,-136){\makebox(0,0)[lb]{\smash{{\SetFigFont{7}{8.4}{\rmdefault}{\mddefault}{\updefault}{\color[rgb]{0,0,0}$\x_C$}%
}}}}
\put(-2324,-136){\makebox(0,0)[lb]{\smash{{\SetFigFont{7}{8.4}{\rmdefault}{\mddefault}{\updefault}{\color[rgb]{0,0,0}$\x_D$}%
}}}}
\put( 76,2939){\makebox(0,0)[lb]{\smash{{\SetFigFont{7}{8.4}{\rmdefault}{\mddefault}{\updefault}{\color[rgb]{0,0,0}$y$}%
}}}}
\put(1951,914){\makebox(0,0)[lb]{\smash{{\SetFigFont{7}{8.4}{\rmdefault}{\mddefault}{\updefault}{\color[rgb]{0,0,0}$x$}%
}}}}
\put(2101,1589){\makebox(0,0)[lb]{\smash{{\SetFigFont{7}{8.4}{\rmdefault}{\mddefault}{\updefault}{\color[rgb]{0,0,0}$\x_B$}%
}}}}
\put(-899,1589){\makebox(0,0)[lb]{\smash{{\SetFigFont{7}{8.4}{\rmdefault}{\mddefault}{\updefault}{\color[rgb]{0,0,0}$\x_A$}%
}}}}
\end{picture}%

%% file: conformal.pstex_t
\begin{picture}(0,0)%
\includegraphics{conformal.pstex}%
\end{picture}%
\setlength{\unitlength}{2486sp}%
\begingroup\makeatletter\ifx\SetFigFont\undefined%
\gdef\SetFigFont#1#2#3#4#5{%
  \reset@font\fontsize{#1}{#2pt}%
  \fontfamily{#3}\fontseries{#4}\fontshape{#5}%
  \selectfont}%
\fi\endgroup%
\begin{picture}(7230,2308)(-1499,-364)
\put(5716,749){\makebox(0,0)[b]{\smash{{\SetFigFont{8}{9.6}{\rmdefault}{\mddefault}{\updefault}{\color[rgb]{0,0,0}$\zeta_C$}%
}}}}
\put(3106,1199){\makebox(0,0)[b]{\smash{{\SetFigFont{8}{9.6}{\rmdefault}{\mddefault}{\updefault}{\color[rgb]{0,0,0}$z = g(\zeta)$}%
}}}}
\put(-269,1559){\makebox(0,0)[b]{\smash{{\SetFigFont{8}{9.6}{\rmdefault}{\mddefault}{\updefault}{\color[rgb]{0,0,0}$z_A$}%
}}}}
\put(5401,1334){\makebox(0,0)[lb]{\smash{{\SetFigFont{8}{9.6}{\rmdefault}{\mddefault}{\updefault}{\color[rgb]{0,0,0}$\Gamma'$}%
}}}}
\put(1531,1559){\makebox(0,0)[b]{\smash{{\SetFigFont{8}{9.6}{\rmdefault}{\mddefault}{\updefault}{\color[rgb]{0,0,0}$z_B$}%
}}}}
\put(-1484, 51){\makebox(0,0)[b]{\smash{{\SetFigFont{8}{9.6}{\rmdefault}{\mddefault}{\updefault}{\color[rgb]{0,0,0}$z_D$}%
}}}}
\put(271, 51){\makebox(0,0)[b]{\smash{{\SetFigFont{8}{9.6}{\rmdefault}{\mddefault}{\updefault}{\color[rgb]{0,0,0}$z_C$}%
}}}}
\put(4726,1761){\makebox(0,0)[b]{\smash{{\SetFigFont{8}{9.6}{\rmdefault}{\mddefault}{\updefault}{\color[rgb]{0,0,0}$\zeta_D$}%
}}}}
\put(4726,-286){\makebox(0,0)[b]{\smash{{\SetFigFont{8}{9.6}{\rmdefault}{\mddefault}{\updefault}{\color[rgb]{0,0,0}$\zeta_B$}%
}}}}
\put(3646,749){\makebox(0,0)[b]{\smash{{\SetFigFont{8}{9.6}{\rmdefault}{\mddefault}{\updefault}{\color[rgb]{0,0,0}$\zeta_A$}%
}}}}
\put(642,1503){\makebox(0,0)[b]{\smash{{\SetFigFont{8}{9.6}{\rmdefault}{\mddefault}{\updefault}{\color[rgb]{0,0,0}$\Gamma$}%
}}}}
\put(474,327){\makebox(0,0)[lb]{\smash{{\SetFigFont{8}{9.6}{\rmdefault}{\mddefault}{\updefault}{\color[rgb]{0,0,0}$\theta$}%
}}}}
\put(1576,614){\makebox(0,0)[b]{\smash{{\SetFigFont{8}{9.6}{\rmdefault}{\mddefault}{\updefault}{\color[rgb]{0,0,0}$z$}%
}}}}
\put(5041,254){\makebox(0,0)[b]{\smash{{\SetFigFont{8}{9.6}{\rmdefault}{\mddefault}{\updefault}{\color[rgb]{0,0,0}$\zeta$}%
}}}}
\end{picture}%